\documentclass[11pt]{amsart}
\usepackage{tocvsec2}
\usepackage{amsmath,amsfonts,latexsym,amssymb,amsthm,mathrsfs,array}

\usepackage[letterpaper]{geometry}
\usepackage[usenames]{color}
\usepackage[colorlinks=true,linkcolor=webgreen,filecolor=webbrown,citecolor=webgreen]{hyperref}

\definecolor{webgreen}{rgb}{0,.5,0}
\definecolor{webbrown}{rgb}{.6,0,0}

\newcommand{\CC}{{\mathbb C}}
\newcommand{\ZZ}{{\mathbb Z}}

\def\gg{{\mathfrak{g}}}
\def\hh{{\mathfrak{h}}}
\def\CC{{\mathbb C}}
\def\ZZ{{\mathbb Z}}

\newtheorem{dfn}{Definition}[section]
\newcommand{\bdfn}{\begin{dfn}\rm}
\newcommand{\edfn}{\end{dfn}}
\newtheorem{thm}[dfn]{Theorem}
\newcommand{\bthm}{\begin{thm}}
\newcommand{\ethm}{\end{thm}}
\newtheorem{exmp}[dfn]{Example}
\newcommand{\bexmp}{\begin{exmp}\rm}
\newcommand{\eexmp}{\end{exmp}}
\newtheorem{lmma}[dfn]{Lemma}                   
\newcommand{\blmma}{\begin{lmma}}                   
\newcommand{\elmma}{\end{lmma}}                   
\newtheorem{ppsn}[dfn]{Proposition}
\newcommand{\bppsn}{\begin{ppsn}}
\newcommand{\eppsn}{\end{ppsn}}
\newtheorem{crlre}[dfn]{Corollary}
\newcommand{\bcrlre}{\begin{crlre}} 
\newcommand{\ecrlre}{\end{crlre}}
\newtheorem{rmk}[dfn]{Remark}
\newcommand{\brmk}{\begin{rmk}\rm} 
\newcommand{\ermk}{\end{rmk}}

\numberwithin{equation}{section}

\title
{Quasi-finite modules
over affine and extended affine Lie algebras}
\author{Souvik Pal}
\address{Department of Mathematics, Indian Institute of Science, CV Raman Road, Bengaluru, Karnataka 560012; and Department of Sciences and Humanities, CHRIST University, Mysore Road, Bangalore 560 074, India.}     
\email{pal.souvik90@gmail.com, souvik.pal@christuniversity.in}      
  
\date{}

\begin{document}



\maketitle
\begin{abstract}
	In this paper, we consider irreducible quasi-finite (or equivalently weakly integrable) modules, with non-trivial action of the core, over the extended affine Lie algebras (EALAs) whose centerless cores are multiloop algebras. The centerless cores of all but one family of EALAs having nullity greater than $1$ are known to admit such multiloop realizations. For any such (untwisted) EALA, we show that the irreducible quasi-finite modules are either integrable with the center of the underlying core acting trivially, or restricted generalized highest weight (GHW) modules. We further prove that in the nullity $2$ case, these irreducible restricted GHW modules turn out to be highest weight type modules, thereby classifying the irreducible quasi-finite modules over all such EALAs. In particular, we obtain the classification of  irreducible quasi-finite modules over toroidal Lie algebras, minimal EALAs and toroidal EALAs of nullity $2$. Along the way, we completely classify the irreducible weakly integrable modules over affine Kac-Moody algebras (studied by Rao--Futorny [\textit{Trans.\ Amer.\ Math.\ Soc.}\ 2009] for non-zero level modules). Our results generalize the well-known work of Chari [\textit{Invent.\ Math.}\ 1986] and Chari--Pressley [\textit{Math.\ Ann.}\ 1986] concerning the classification of irreducible integrable modules over (nullity $1$) affine Kac--Moody algebras.    
\end{abstract}

\settocdepth{section}
\tableofcontents
\section{Introduction}
The affine Kac--Moody algebras are an important family of infinite-dimensional Lie algebras, for which there exists a well-developed structure theory and representation theory that effectively mirrors the classical theory of finite-dimensional simple Lie algebras \cite{K}. Most of the applications of affine Kac--Moody algebras arise from their representation theory, which surprisingly turns out to be a powerful tool for the investigation of several apparently disconnected fields, like vertex algebras, modular forms and the Geometric Langlands program, to name a few \cite{K,MV}.

The construction of an affine Kac--Moody algebra proceeds in two  steps, which we recall here. For simplicity, we shall henceforth work over the complex numbers $\mathbb{C}$, although all our results also hold over any algebraically closed field of characteristic $0$. Let $\mathfrak{g}$ be a finite-dimensional simple Lie algebra with a Cartan subalgebra $\mathfrak{h}$. Then the (untwisted) affine Kac--Moody algebra is given by
\begin{align*}
\widehat{L}(\mathfrak{g}) = \mathfrak{g} \otimes \mathbb{C}[t,t^{-1}] \oplus \mathbb{C}K \oplus \mathbb{C}d,
\end{align*}  
where $\mathfrak{g} \otimes \mathbb{C}[t,t^{-1}] \oplus \mathbb{C}K$ is a central extension (in fact, the universal central extension) of $\mathfrak{g} \otimes \mathbb{C}[t,t^{-1}]$ and $d$ is the degree derivation $t(d/dt)$ acting on $\mathbb{C}[t,t^{-1}]$, which just tracks the $\mathbb{Z}$-grading of $\widehat{L}(\mathfrak{g})$.

Extended affine Lie algebras (EALAs) are a rich class of Lie algebras which initially appeared in the context of elliptic singularities \cite{S} and in the framework of Lie algebras related to quantum field gauge theory \cite{KT}. These are natural higher-dimensional generalizations of the finite-dimensional simple Lie algebras and affine Kac--Moody algebras. For instance, replacing $A_1 = \mathbb{C}[t,t^{-1}]$ by $A_{n+1} = \mathbb{C}[t_0^{\pm 1}, \ldots, t_n^{\pm 1}]$ and then  choosing appropriate central extensions and derivations gives rise to concrete examples of EALAs. However, there are also other examples of EALAs which can be formed by considering a more general family of Lie algebras known as the \textit{Lie torus}, instead of $\mathfrak{g} \otimes A_{n+1}$ (see \S \ref{SS2.1}). In fact, it was shown by Neher \cite{Ne} that every EALA can be constructed from a Lie torus by taking a central extension and then adding a suitable space of derivations.
 
Each EALA has an invariant non-negative integer attached to it, called its \textit{nullity}. The EALAs of nullity $0$ are simply the finite-dimensional simple Lie algebras, whereas the nullity $1$ EALAs are precisely the affine Kac--Moody algebras. It is also worth mentioning  that the EALAs of nullity $2$ have natural connections with double affine Hecke algebras (DAHAs) via double affine Weyl groups. Similar to the nullity $0$ and nullity $1$ cases, one can associate a root system to any EALA, which can be partitioned into the sets of isotropic and non-isotropic roots (see \S \ref{SS3.1}). The subalgebra generated by the non-isotropic root vectors forms an ideal and is known as the \textit{core of an EALA}, which is always a Lie torus.

Let $E_{n+1}$ be an EALA of nullity $n+1$ ($\geqslant 1$) with core $(E_{n+1})_c$. Due to remarkable breakthroughs in \cite{ABFP,Ne}, it is now evident that except for a particular class, namely when the root system of the underlying centerless core $(E_{n+1})_{cc}$ of $E_{n+1}$ is not of type $A$, or more precisely, when $(E_{n+1})_{cc}$ does not satisfy the \textit{fgc condition} (see \S \ref{SS2.2} and \S \ref{SS3.1}), $E_{n+1}$ can be always realized as an extension of a twisted multiloop algebra (or a twisted loop algebra for $n =0$). In this particular set-up, if the multiloop algebra $L(\mathfrak{g}, \underline{\sigma})$ (see \S \ref{SS2.3}) is untwisted (i.e. $L(\mathfrak{g}, \underline{\sigma}) = \mathfrak{g} \otimes A_{n+1}$), then we shall refer to the corresponding EALAs as (untwisted) \textit{fgc EALAs} (see \S \ref{SS3.3}). Consequently 
\begin{align*}
E_{n+1} = L(\mathfrak{g}) \oplus \mathcal{C}_{n+1} \oplus \mathcal{D}_{n+1}, \ \text{where} \ L(\mathfrak{g}) = \mathfrak{g} \otimes A_{n+1}. 
\end{align*}  
Here $(E_{n+1})_c = L(\mathfrak{g}) \oplus \mathcal{C}_{n+1}$ is a central extension of $L(\mathfrak{g})$ and $\mathcal{D}_{n+1}$ is a $\mathbb{Z}^{n+1}$-graded subalgebra of $\mathcal{S}_{n+1}$ containing the \textit{abelian ad-diagonalizable} subalgebra $D = \mathrm{span} \{d_0, \ldots, d_n \}$ (see \S \ref{SS3.3}). 

Before the advent of EALAs, the \textit{toroidal Lie algebra} (see \S \ref{SS4.2}) was also introduced out of an attempt to generalize the theory of affine Lie algebras in the higher-dimensional framework. This Lie algebra is formed by first taking the universal central extension of $L(\mathfrak{g})$ and then adding $D$, i.e.
\begin{align*}
\tau_{n+1}^{tor} = L(\mathfrak{g}) \oplus \mathcal{Z}_{n+1} \oplus D, \ \text{where} \ \mathrm{dim}\mathcal{Z}_{n+1} = \infty \ \forall \ n \geqslant 1.
\end{align*}      
The structure theory of EALAs has been extensively studied over the last $25$ years (see \cite{AABGP,ABFP,CNPY,Ne} and the references therein). But their representation theory is not so well understood, even in the lower nullity cases. The most important and widely studied class of representations are probably the ones having finite-dimensional weight space decompositions, thereby admitting characters. These modules are usually referred to as \textit{quasi-finite} (or Harish-Chandra) modules in the literature.

In the nullity $0$ case, the irreducible quasi-finite modules over $\mathfrak{g}$ were classified by O. Mathieu \cite{Ma}. In the present paper, we consider the analogous question for $\tau_{n+1}^{tor}$ and fgc EALAs of nullity at least $1$. Unlike the nullity $0$ EALAs, $E_{n+1}$ (as well as $\tau_{n+1}^{tor}$) is $\mathbb{Z}^{n+1}$-graded and so a natural class of quasi-finite modules over $E_{n+1}$ (respectively $\tau_{n+1}^{tor}$) are the weight modules having finite-dimensional graded components (or equivalently having finite-dimensional $D$-weight spaces). 

The irreducible \textit{integrable} modules over $\widehat{L}(\mathfrak{g})$ and $\tau_{n+1}^{tor}$ were classified in \cite{C,CP} and \cite{E1,E2} respectively, whereas the irreducible \textit{weakly integrable} $\widehat{L}(\mathfrak{g})$-modules of \textit{non-zero level} were studied in \cite{EF} (see Section \ref{S7} for definitions). In the higher nullity set-up, the irreducible \textit{integrable}  modules over $E_{n+1}$ (with non-trivial $(E_{n+1})_c$-action) were classified in some isolated cases, like toroidal EALAs \cite{CLT,BES,TB}, minimal EALAs \cite{E1} and Hamiltonian EALAs \cite{E4}, obtained by making specific choices of $\mathcal{C}_{n+1}$ and $\mathcal{D}_{n+1}$ (see Section \ref{S4}). Vertex operator realizations of several quasi-finite modules over $\tau_{n+1}^{tor}$ and toroidal EALAs were also provided in \cite{EMY} and \cite{YB,FLT} respectively. However, the classification problem of such irreducible quasi-finite (or equivalently \textit{weakly integrable}) modules over $E_{n+1}$ is still unresolved, even in the nullity $1$ case (see Remark \ref{R8.4}).

In this paper, we attempt to address the aforesaid classification problems in a unified manner, especially in the nullity $1$ and nullity $2$ cases. More specifically, we establish the following. 
\begin{itemize}
\item Classify the irreducible quasi-finite modules over the affine Kac--Moody algebra $\widehat{L}(\mathfrak{g})$.
\item Classify the irreducible quasi-finite modules (with non-trivial core action) over $\tau_2^{tor}$ as well as over the nullity $2$ fgc EALAs by pursuing a unified approach.
\item Broadly describe the irreducible quasi-finite modules (with non-trivial action of the core) over $\tau_{n+1}^{tor}$ and $E_{n+1}$ for $n>1$.	
\end{itemize}
The precise statements related to the above classification problems can be found in Section \ref{S11}, which generalize the results of Chari--Pressley \cite{C,CP} related to the classification of irreducible \textit{integrable} modules over $\widehat{L}(\mathfrak{g})$ and also the work of Rao--Futorny \cite{EF}, where the irreducible \textit{weakly integrable} $\widehat{L}(\mathfrak{g})$-modules of \textit{non-zero level} were shown to be highest weight modules.  
 
We finally remark that the representation theory of the above three classes of EALAs, namely the affine Kac--Moody algebras, fgc EALAs of nullity $2$ and fgc EALAs of nullity greater than $2$ are vastly different from each other (for example, see \cite[Remark 5.13]{CLT} and \cite{BES}). The main obstacle that we face while working with $E_{n+1}$ ($n \geqslant 1$) is the highly arbitrary nature of $\mathcal{D}_{n+1}$, due to which we cannot effectively utilize the representation theory of $\mathcal{D}_{n+1}$ to tackle our problem.

We hope that our results will contribute towards solving the more general problem of classifying irreducible quasi-finite modules over EALAs as well as developing the overall representation theory of EALAs. In a recent paper \cite{Ku}, irreducible highest weight modules over minimal EALAs (see \S \ref{SS4.4}) were studied in connection to modular representation theory, while addressing some questions raised by Lusztig. Another promising direction is the classification of irreducible weakly integrable modules over affine Lie superalgebras (see \cite{KW,EF} for partial results). 

\smallskip
   	
\noindent \textbf{Organization of the paper.} In Section \ref{S2}, we record some results related to a Lie torus which are relevant to our paper. In Section \ref{S3}, we recall Neher's general construction of an EALA from its centerless core (see Theorem \ref{T3.1}) and also present the explicit realization of an (untwisted) EALA whose centerless core satisfies fgc condition (see Corollary \ref{C3.3}). In Lemma \ref{L3.6}, we prove an important result related to the core and center of an  fgc EALA. In Section \ref{S4}, we provide some concrete examples of fgc EALAs after recalling the notions of toroidal and full toroidal Lie algebras. In Section \ref{S5}, we define the Weyl group of an fgc EALA along with its roots and co-roots, while in Section \ref{S6}, we discuss automorphisms of fgc EALAs.              

In Section \ref{S7}, we prove an important lemma regarding the action of the central elements of $E_{n+1}$ (see Lemma \ref{L7.4}) and also show that the category of quasi-finite modules coincides with the category of weakly integrable modules over $E_{n+1}$ (see Lemma \ref{L7.8}). In Section \ref{S8}, we completely classify the irreducible quasi-finite $\widehat{L}(\mathfrak{g})$-modules (see Theorem \ref{T8.5}). In Section \ref{S9}, we show that any irreducible quasi-finite module over $E_{n+1}$ (or $\tau_{n+1}^{tor}$) with non-trivial action of $(E_{n+1})_c$, is either a \textit{restricted generalized highest weight} (GHW) module or an \textit{integrable} highest weight module over $L(\mathfrak{g})  \oplus \mathcal{D}_{n+1}$, which is induced from an irreducible \textit{uniformly bounded} module over $\mathfrak{h} \otimes A_{n+1} \oplus \mathcal{D}_{n+1}$ (see Proposition \ref{P9.5}, Remark \ref{R9.6} and Proposition \ref{P9.9}). In Section \ref{S10}, we prove that the irreducible restricted GHW modules over $E_2$ are in fact \textit{highest weight type modules} 
(see Theorem \ref{T10.15}), thereby classifying all the irreducible quasi-finite modules over $E_2$ (and $\tau_2^{tor}$) with non-trivial core action. In Section \ref{S11}, we record the main results of our paper. As an application of Theorem \ref{T11.3}, we obtain the classification of irreducible quasi-finite modules, with non-trivial core action, over minimal EALAs and toroidal EALAs of nullity $2$ as well as over $\tau_2^{tor}$ (see Remark \ref{R11.4}), whereas for $n>1$, we broadly characterize these irreducible quasi-finite modules over $E_{n+1}$ and $\tau_{n+1}^{tor}$ (see Theorem \ref{T11.5} and Remark \ref{R11.6}). 

\smallskip
	    
\noindent \textbf{Acknowlegements.} The author would like to thank Prof. Eswara Rao and Dr. G. Krishna Teja for some discussions. He also thanks Prof. Apoorva Khare for going through the whole paper and for his helpful suggestions regarding its exposition. 
Most of this work was done when the author was a Research Associate at the Theoretical Statistics and Mathematics Unit, Indian Statistical Institute, Bangalore, India. The author gratefully acknowledges the financial support and the excellent working conditions provided by the institute.                              

\smallskip

\noindent \textbf{Notations.} Throughout this paper, all the vector spaces, algebras and tensor products are over the field of complex numbers $\mathbb{C}$. We shall denote the set of integers, natural numbers, non-negative integers and non-zero complex numbers by $\mathbb{Z}$, $\mathbb{N}$, $\mathbb{Z_{+}}$ and $\mathbb{C}^{\times}$ respectively. For any Lie algebra $L$, its universal enveloping algebra will be denoted by $U(L)$  and its center will be denoted by $Z(L)$. 

\section{Lie torus and its multiloop realization} \label{S2}
In this section, we recall the notion of a Lie torus satisfying the fgc condition and then provide its multiloop realization, which will be used throughout this paper. We conclude by briefly discussing about the derivations arising from a Lie torus. The interested reader is referred to \cite{ABFP,N1,Ne,EN} for more details.
 
\subsection{Definition and examples} \label{SS2.1}
Let $\Delta$ be a finite irreducible (but not necessarily reduced) root system including $0$ and $\Lambda$ be any abelian group. Put $\Delta_{\text{ind}} = \{0 \} \cup \{\alpha \in \Delta \ | \ \alpha/2 \notin \Delta \}$. Then a \textit{Lie torus of type $(\Delta, \Lambda)$} is a Lie algebra $L$ satisfying the following properties. 

\noindent (LT1) (a) $L$ is $\Lambda$-graded as well as $Q(\Delta)$-graded, where $Q(\Delta) := \mathbb{Z}\Delta$ is the root lattice of $\Delta$. Then 

$\,\,\,\,\,\,\,\,\,\,\,\,\,\,\,\ $ $L = \bigoplus_{\alpha \in Q(\Delta), \lambda \in \Lambda} L_{\alpha}^{\lambda}$, with $[L_{\alpha}^{\lambda}, L_{\beta}^{\mu}] \subseteq L_{\alpha + \beta}^{\lambda + \mu}$. Set $L_{\alpha} = \bigoplus_{\lambda \in \Lambda} L_{\alpha}^{\lambda}$ and $L^{\lambda} = \bigoplus_{\alpha \in Q(\Delta)} L_{\alpha}^{\lambda}$. 

$\,\,\,\,\,\,\,\ $ (b) $\text{supp}_{Q(\Delta)} L = \{\alpha \in Q(\Delta) \ | \ L_{\alpha} \neq (0) \} = \Delta$, which gives $L = \bigoplus_{\alpha \in \Delta} L_{\alpha}$. 

\noindent (LT2) (a) If $\alpha \neq 0$ and $L_{\alpha}^{\lambda} \neq (0)$, then there exist some $e_{\alpha}^{\lambda} \in L_{\alpha}^{\lambda}$ and $f_{\alpha}^{\lambda} \in L_{-\alpha}^{-\lambda}$ such that $L_{\alpha}^{\lambda} =
 \mathbb{C}e_{\alpha}^{\lambda}$, 
 
 $\,\,\,\,\,\,\,\,\,\,\,\,\,\,\,\ $ $L_{-\alpha}^{-\lambda} = \mathbb{C}f_{\alpha}^{\lambda}$ and $[[e_{\alpha}^{\lambda}, f_{\alpha}^{\lambda}], x_{\beta}] = \beta(\alpha^{\vee})x_{\beta} \ \forall \ \beta \in \Delta, \ x_{\beta} \in L_{\beta}$, where $\alpha^{\vee}$ denotes the co-root
 
$\,\,\,\,\,\,\,\,\,\,\,\,\,\,\,\ $  corresponding to $\alpha$ in the usual sense.    

$\,\,\,\,\,\,\,\ $ (b) $L_{\alpha}^0 \neq (0) \ \forall \ 0 \neq \alpha \in \Delta_{\text{ind}}$.

\noindent (LT3) $L$ is generated by $\bigcup_{0 \neq \alpha \in \Delta} L_{\alpha}$ as a Lie algebra. 

\noindent (LT4) As an abelian group, $\Lambda$ is generated by supp$_{\Lambda}L = \{\lambda \in \Lambda \ | L^{\lambda} \neq (0) \}$. 

\smallskip

\noindent The rank of $\Lambda$ is said to be the \textit{nullity} of $L$. An important example of a Lie torus of type $(\Delta_{\mathfrak{g}}, \mathbb{Z}^{n+1})$ is the (untwisted) multiloop algebra $\mathfrak{g} \otimes \mathbb{C}[t_0^{\pm 1}, \ldots, t_n^{\pm 1}]$, where $\mathfrak{g}$ is a finite-dimensional simple Lie algebra and $\mathbb{C}[t_0^{\pm 1}, \ldots, t_n^{\pm 1}]$ is the Laurent polynomial algebra (see \S \ref{SS2.3}). Another notable example of a Lie torus is $\mathfrak{sl}_{n+1}(\mathbb{C}_{\textbf{q}})$, where $\mathbb{C}_{\textbf{q}}$ is a quantum torus, which has been studied in \cite{BGK}. 

\subsection{Properties of Lie tori} \label{SS2.2}
Let $L$ be a Lie torus of type $(\Delta, \Lambda)$ of nullity $n+1$.
Then:
\smallskip

\noindent (P1) $L$ has a non-zero invariant symmetric bilinear form $\langle \cdot, \cdot \rangle$, which is \textit{$\Lambda$-graded} in the sense that $\langle L^{\lambda}, L^{\mu} \rangle = 0$ if $\lambda + \mu \neq 0$ \cite{Y3}.  \\
(P2) $L$ is a perfect Lie algebra and the universal central extension of $L$ is again perfect. \\
(P3) The \textit{centroid} of $L$, given by $\mathrm{Cent}(L) = \{T \in \mathrm{End}L \ | \ [T, \mathrm{ad}x] = 0 \ \forall \ x \in L \} = \bigoplus_{\lambda \in \Lambda} (\mathrm{Cent}(L))^{\lambda}$ is $\Lambda$-graded with $\mathrm{dim} (\mathrm{Cent}(L))^{\lambda} \leqslant 1$, where $(\mathrm{Cent}(L))^{\lambda}$ consists of endomorphisms of degree $\lambda$ with respect to the $\Lambda$-grading of $L$. Clearly $L$ can be thought of as a module over $\mathrm{Cent}(L)$ under the natural action. If 
we now put $\Gamma = \{\lambda \in \Lambda \ | \ (\mathrm{Cent}(L))^{\lambda} \neq (0) \}$, then $\Gamma$ is a subgroup of $\Lambda$ and is known as the \textit{central grading group}. It can be further shown that $\mathrm{Cent}(L)$  is isomorphic to the group ring $\mathbb{C}[\Gamma]$ and hence is a Laurent polynomial algebra in $k$ variables with $0 \leqslant k \leqslant n+1$. This allows us to express $\mathrm{Cent}(L) = \bigoplus_{\lambda \in \Gamma} \mathbb{C}\chi^{\lambda}$, where the $\chi^{\lambda}$'s satisfy the multiplication rule $\chi^{\lambda_1} \chi^{\lambda_2} = \chi^{\lambda_1 + \lambda_2}$ and act on $L$ as endomorphisms of $\Lambda$-degree $\lambda$. \\
(P4) If $\Delta$ is not of type $A_l$ for any $l \in \mathbb{N}$, then it was shown in \cite{N1} that $L$ has finite rank as a $\text{Cent}(L)$-module or equivalently $k = n+1$, where $k$ is as given in (P3). In this case, we say that $L$ satisfies the \textit{fgc condition}, i.e. $L$ is finitely generated over its centroid. Lie tori of type $A_l$ were classified in \cite{BGK,BGKN,Y1}. \\
(P5) If $L$ does not satisfy the fgc condition, then it follows from (P4) and \cite{NY} that $L \cong \mathfrak{sl}_l(\mathbb{C}_{\textbf{q}})$, where $\mathbb{C}_{\textbf{q}}$ is a quantum torus in $n+1$ variables and $\textbf{q} = (q_{ij})$ is an $(n+1) \times (n+1)$ quantum matrix with at least one $q_{ij}$ not a root of unity.    

\subsection{Multiloop algebras} \label{SS2.3} 
Suppose $\mathfrak{g}$ is a finite-dimensional simple Lie algebra with the usual triangular decomposition $\mathfrak{g} = \mathfrak{n}^{-} \oplus \mathfrak{h} \oplus \mathfrak{n}^{+}$ relative to a Cartan subalgebra $\mathfrak{h}$. It is well-known that $\mathfrak{g}$ is endowed with a symmetric, non-degenerate and invariant bilinear from $(\cdot| \cdot)$, which is unique up to scalars. Let $A_{n+1} = \CC[t_0^{\pm 1}, \ldots, t_n^{\pm 1}]$ ($n \in \mathbb{Z}_{+}$) be the algebra of Laurent polynomials in $n+1$ variables.
Consider the (untwisted) multiloop algebra given by
\[ {
	L(\mathfrak{g}) =\mathfrak{g} \otimes A_{n+1} , \qquad
	[x \otimes f, y \otimes g] = [x,y] \otimes fg \ \forall \ x,y \in \mathfrak{g} \ \text{and} \ f, g \in A_{n+1}.
}
\]

\noindent Fix any $n\in\mathbb{Z}_{+}$ and suppose we have $n+1$ commuting automorphisms of $\mathfrak{g}$  given by $\sigma_0, \ldots, \sigma _n$ with finite orders $m_0, \ldots, m_n$ respectively. We choose an $l$-th primitive root of unity $\xi_l$ for all positive integers in such a way that $\xi_{ml}^m = \xi_l$ for $m, l \in \mathbb{N}$. Put
\begin{align*} 
\underline{\sigma} = (\sigma_0, \ldots, \sigma_n),\,\,\,\,\,\,\,\ \Lambda = \mathbb{Z}^{n+1}, \,\,\,\,\,\,\,\,\
\Gamma = m_0 \mathbb{Z} \oplus \ldots \oplus m_n\mathbb{Z},\,\,\,\,\,\,\,\
G = \Lambda/ \Gamma.
\end{align*} 
Thus we have a natural map \ $\Lambda \longrightarrow G$ $(\cong \mathbb{Z}/m_0\mathbb{Z} \times \ldots \times\mathbb{Z}/m_n\mathbb{Z})$ 
\begin{align}\label{Natural}
(k_0, \ldots, k_n) = \underline{k} \quad \mapsto \quad \overline{k}=(\overline{k_0}, \ldots ,\overline{k_n})
\end{align} 
For $0 \leqslant i \leqslant n$, let  $\omega_i$ denote an $m_i-$th primitive root of unity. Then we obtain an eigenspace decomposition of $\mathfrak{g}$ given by
\begin{align}\label{Eigenspace}
\mathfrak{g} = \bigoplus _{\overline{k} \in G} \mathfrak{g}_{\overline{k}} \,\,\,\ \text{where} \,\,\,\ \mathfrak{g}_{\overline{k}} := \{x \in \mathfrak{g} \ | \ \sigma_i x = \omega_i^{k_i}x,\ 0 \leqslant i \leqslant n \}.
\end{align}
The subalgebra $\mathfrak{g}_{\overline{0}}$ is a reductive Lie algebra, even with the possibility of being zero. Finally, define
\begin{align*} 
L(\mathfrak{g}, \underline{\sigma}) := \bigoplus_{\underline{k} \in \Lambda} \mathfrak{g}_{\overline{k}} \otimes \mathbb{C}t^{\underline{k}},
\end{align*}
which is known as the \textit{twisted multiloop algebra} associated to $\mathfrak{g}$ and $\underline{\sigma}$.

\brmk
In general, $L(\mathfrak{g}, \underline{\sigma})$ is not necessarily a Lie torus. Nevertheless, a careful choice of the automorphisms $\underline{\sigma}$ of $\mathfrak{g}$ does indeed give rise to a twisted multiloop algebra which is also a Lie torus (see \cite[Proposition 3.2.5]{ABFP} for the precise conditions). Such a choice ensures that $\mathfrak{g}_{\overline{0}}$ is a (non-zero) simple Lie algebra and in this case, $L(\mathfrak{g},\sigma)$ coincides with the twisted loop algebra used in the construction of twisted affine Kac--Moody algebras (for $n = 0$).
\ermk

The root space decomposition of $\gg$ relative to $\hh$ is given by
\begin{align*}
\gg = \bigoplus_{\alpha \in \hh^*}\mathfrak{g}_{\alpha} \,\,\,\  \text{where} \,\,\,\  \mathfrak{g}_{\alpha} = \{ x \in \gg \ | \ [h,x]= \alpha(h)x \ \forall \ h \in \hh \}.     
\end{align*}
Set $\Delta_{\gg} := \{ \alpha \in \hh^*\ |\ \mathfrak{g}_{\alpha} \neq (0)\}$. Then $\Delta_{\mathfrak{g}}^{\times} = \Delta_{\gg} \backslash \{0\}$ is an irreducible reduced finite root system.
Let $\Delta_{\mathfrak{g}}^{+} \subseteq \Delta_{\mathfrak{g}}^{\times}$ denote the set of all positive roots of $\mathfrak{g}$ with respect to $\mathfrak{h}$ and $Q$ be the corresponding root lattice. Then the Lie torus $L(\mathfrak{g})$ is $\mathbb{Z}^{n+1}$-graded as well as $Q$-graded and is of type $(\Delta_{\mathfrak{g}}, \mathbb{Z}^{n+1})$. We shall denote the set of all dominant integral weights of $\mathfrak{g}$ by
\begin{align*}
P^{+}_{\mathfrak{g}} = \{ \lambda\ \in \mathfrak{h}^*\ | \ (\lambda|\alpha) \in \mathbb{Z}_{+}\ \forall \ \alpha \in \Delta_{\mathfrak{g}}^{+} \}.
\end{align*}

\bdfn \label{D2.2}
$\mu \in \mathfrak{h}^*$ is said to be a minimal (or miniscule) weight of $\mathfrak{g}$ if $\mu \in P_{\mathfrak{g}}^{+}$ and $\mu(\alpha^{\vee}) \in \{-1,0,1 \} \ \forall \ \alpha \in \Delta_{\mathfrak{g}}^{\times}$, where $\alpha^{\vee} \in \mathfrak{h}$ denotes the co-root corresponding to the root $\alpha$.
\edfn

\subsection{Centroidal derivations of Lie tori(\cite{N1,Ne})} \label{SS2.4} 
Let $L$ be a centerless Lie torus of type $(\Delta, \Lambda)$ with central grading group $\Gamma$ (see (P3) of \S \ref{SS2.2}). Now each $\theta \in \mathrm{Hom}_{\mathbb{Z}}(\Lambda, \mathbb{C})$ gives rise to a \textit{degree derivation} $\partial_{\theta}$ of $L$, defined by $\partial_{\theta}(x^{\lambda}) = \theta(\lambda)x^{\lambda}$ for $x^{\lambda} \in L^{\lambda}$. Clearly $D := \{\partial_{\theta} \ | \ \theta \in \mathrm{Hom}_{\mathbb{Z}}(\Lambda,\mathbb{C}) \}$ gives rise to the $\Lambda$-grading of $L$ given by $L^{\lambda} := \{x \in L \ | \ \partial_{\theta}(x) = \theta(\lambda)x \ \forall \ \partial_{\theta} \in D \}$. Let us denote the algebra of derivations of $L$ by $\mathrm{Der}(L) := \{\delta \in \mathrm{End}_{\mathbb{C}}(L) \ | \ \delta([x,y]) = [x,\delta(y)] + [\delta(x),y] \ \forall \ x, y \in L \}$. Now if $T \in \mathrm{Cent}(L)$, then it is easy to see that $Td \in \mathrm{Der}(L)$ for any $d \in \mathrm{Der}(L)$. As a result, 
\begin{align*}
\mathrm{CDer}(L) := \mathrm{Cent}(L)D = \bigoplus_{\lambda \in \Gamma} \chi^{\lambda}D
\end{align*}
is a $\Gamma$-graded subalgebra of $\mathrm{Der}(L)$, which is known as the algebra of \textit{centroidal derivations} of $L$. The bilinear operation on $\mathrm{CDer}(L)$ is then given by
\begin{align*} 
[\chi^{\lambda_1} \partial_{\theta_1}, \chi^{\lambda_2} \partial_{\theta_2}] = \chi^{\lambda_1 + \lambda_2}(\theta_1(\lambda_2)\partial_{\theta_2} - \theta_2(\lambda_1) \partial_{\theta_1}), 
\end{align*}       
whence $\mathrm{CDer}(L)$ is a $\Gamma$-graded subalgebra of $\mathrm{Der}(L)$ with an ad-diagonalizable subalgebra $D$.

The $\Gamma$-graded subalgebra of \textit{skew-centroidal derivations}  with respect to $<\cdot, \cdot>$ is defined as
\begin{align*}
\mathrm{SCDer}(L) := \{d \in \mathrm{CDer}(L) \ | \ \langle d(l_1), l_2 \rangle + \langle l_1, d(l_2) \rangle = 0 \ \forall \ l_1, l_2 \in L \} \\ 
= \bigoplus_{\lambda \in \Gamma}(\mathrm{SCDer}(L))^{\lambda}, \ \text{where} \ (\mathrm{SCDer}(L))^{\lambda} = \chi^{\lambda} \{\partial_{\theta} \in D \ | \ \theta(\lambda) = 0 \}. 
\end{align*} 
Then it is evident that $\mathrm{SCDer}(L)^{0} = D$ and $[(\mathrm{SCDer}(L))^{\lambda}, (\mathrm{SCDer}(L))^{-\lambda}] = (0)$, due to which $\mathrm{SCDer}(L)$ can be written as a semi-direct product of $D$ and $\bigoplus_{0 \neq \lambda \in \Gamma} (\mathrm{SCDer}(L))^{\lambda}$.  
      
\section{Fgc extended affine Lie algebras} \label{S3}
In this section, we recall the general definition of extended affine Lie algebras and then provide the explicit realization of fgc extended affine Lie algebras following \cite{CNPY,N1,Ne,EN}.
At the end of the section, we prove an important lemma which will play a vital role throughout this paper.

\subsection{Axiomatic definition of general EALAs} \label{SS3.1} 
 An extended affine Lie algebra (EALA for short) consists of a triplet $(E,H,\langle \cdot|\cdot \rangle)$, where $E$ is a Lie algebra, $H$ is a subalgebra of $E$ and $\langle\cdot|\cdot\rangle$ is a symmetric, non-degenerate, invariant bilinear form on $E$ (in the sense that $\langle [x,y]|z \rangle = \langle x|[y,z] \rangle \ \forall \ x,y,z \in E$) satisfying the following axioms. 
 
\smallskip 

\noindent (EA1) $H$ is a non-trivial finite-dimensional ad-diagonalizable and self-centralizing subalgebra of $E$. \\
In order to describe the remaining axioms, we need to introduce some notations. \\
By (EA1), we have
\begin{align*}
E = \bigoplus_{\alpha \in H^*} E_{\alpha}, \ \text{where} \ E_{\alpha} := \{x \in E \ | \ [h,x] = \alpha(h)x \ \forall \ h \in H \}.
\end{align*}
Let $R := \{\alpha \in H^* | \ E_{\alpha} \neq (0) \}$ be the roots of $E$ with respect to $H (= E_0)$. Note that $0 \in R$. Also $\langle\cdot|\cdot \rangle$ remains non-degenerate when restricted to $H$ and so we can transfer this form to $H^*$. Set
\begin{align*}
R^{\times}:= \{\alpha \in R \ | \ \langle \alpha|\alpha \rangle \neq 0 \}, \ R^0 := \{\alpha \in R \ | \ \langle \alpha|\alpha \rangle = 0 \}.
\end{align*}
The elements of $R^{\times}$ (respectively $R^0$) are called non-isotropic (respectively isotropic) roots. Then we have $R = R^{\times} \sqcup R^0$. The \textit{core} of $E$, which we shall denote by $E_c$, is defined as the subalgebra generated by $\bigcup_{\alpha \in R^{\times}} E_{\alpha}$. It is worth mentioning here that $E_c$ also forms an ideal of $E$. \\
(EA2) For each $\alpha \in R^{\times}$ and $x_{\alpha} \in E_{\alpha}$, ad$x_{\alpha}$ acts locally nilpotently on $E$. \\
(EA3) $R$ is a discrete subset of $H^*$. \\
(EA4) $R$ is an irreducible root system, in the sense that \\
(a) for any decomposition $R^{\times} = R_1 \cup R_2$ satisfying $\langle R_1|R_2 \rangle = (0)$, either $R_1 = (0)$ and $R_2 = (0)$; \\
(b) if $\beta \in R^0$, then there exists $\alpha \in R^{\times}$ such that $\beta + \alpha \in R$. \\
(EA5) The $\mathbb{Z}$-subgroup $\Lambda$ generated by $R^0$ in $(H^*,+)$ is a free abelian group of finite rank. 

It can be deduced from \cite{N1,Ne} that there exists a finite irreducible root system $\Delta$, an embedding $\Delta_{\text{ind}} \hookrightarrow R$ and a family $\{\Lambda_{\alpha} \}_{\alpha \in \Delta} \subseteq \Lambda$ such that $R = \bigcup_{\alpha \in \Delta} (\alpha + \Lambda_{\alpha})$ and consequently $E_c$ is a Lie torus of type $(\Delta, \Lambda)$. Then the centerless core of $E_c$, given by $E_{cc} = E_c/Z(E_c)$, is again a (centerless) Lie torus, having the same type and nullity as $E_c$ \cite{Y2}.

The rank of $\Lambda$ is known as the \textit{nullity} of $E$. Any EALA of nullity $0$ is a finite-dimensional simple Lie algebra \cite{EN}, whereas the nullity $1$ EALAs are precisely the affine Kac--Moody algebras \cite{ABGP}.   

\subsection{General construction of EALAs from Lie tori} \label{SS3.2}
Let $L$ be a centerless Lie torus of type $(\Delta, \Lambda)$ and nullity $n+1$, where $n \geqslant 0$. Suppose $\langle \cdot,\cdot \rangle$ is a $\Lambda$-graded invariant symmetric non-degenerate bilinear form on $L$ and let $\Gamma$ be the central grading group. Let $\mathcal{D}_{n+1} := \mathcal{D} = \bigoplus_{\lambda \in \Gamma} \mathcal{D}_{\lambda}$ be a $\Gamma$-graded subalgebra of $\mathrm{SCDer}(L)$ such that the canonical evaluation map ev: $\Lambda \longrightarrow \mathcal{D}_0^*$ given by (ev($\lambda$))($\partial_{\theta}$) = $\theta(\lambda)$ is injective. This is equivalent to the fact that $\mathcal{D}_0$ induces the $\Lambda$-grading of $L$ (with $L^{\lambda} = \{x \in L \ | \ \partial_{\theta}(x) = \theta(\lambda)x \ \forall \ \partial_{\theta} \in \mathcal{D}_0 \} \ \forall \ \lambda \in \Lambda$) and in particular implies that $\mathcal{D}_0 \neq (0)$. If $\mathcal{C}_{n+1} := \mathcal{C} = \bigoplus_{\lambda \in \Gamma} \mathcal{C}^{\lambda}$ denotes the graded dual space of $\mathcal{D}$, then we can consider $\mathcal{C}$ as a $\Gamma$-graded vector space with the $\lambda$-th component $\mathcal{C}_{\lambda} = (\mathcal{D}_{-\lambda})^*$, which reveals that $\mathrm{dim} \ \mathcal{C}_0 = \mathrm{dim} \ \mathcal{D}_0 \geqslant 1$. Also $\mathcal{D}$ acts on $L$ via derivations of $L$, which we shall denote by $d(l)$. Then $\sigma_{\mathcal{D}}(l_1,l_2)(d) = \langle d(l_1),l_2 \rangle$ is a central $2$-cocycle for $L$ with values in $\mathcal{C}$. Clearly $\mathcal{C}$ remains invariant under the (graded) contragredient action of $\mathcal{D}$ on $\mathcal{D}^*$, written as $d \cdot f$. The next ingredient needed for our construction is an affine $2$-cocycle $\phi : \mathcal{D} \times \mathcal{D} \longrightarrow \mathcal{C}$ satisfying
\begin{align} \label{(3.1)}
\phi(d_0,d) = 0 \ \forall \ d_0 \in \mathcal{D}_0, \ d \in \mathcal{D} \,\,\,\ \text{and} \,\,\,\ \phi(d_1,d_2)(d_3) = \phi(d_2,d_3)(d_1) \ \forall \ d_1, d_2, d_3 \in \mathcal{D}. 
\end{align} 
Finally, set
\begin{align} \label{(3.2)}
E := E(L, \mathcal{D}, \phi) = L \oplus \mathcal{C} \oplus \mathcal{D}.
\end{align}
It was shown in \cite{CNPY} that $E$ forms a Lie algebra with respect to the following bracket operation:
\begin{align*}
[l_1 + f_1 + d_1, l_2 + f_2 + d_2] := \\ ([l_1,l_2]_{L} + d_1(l_2) - d_2(l_1)) + 
(\sigma_{\mathcal{D}}(l_1,l_2) + d_1 \cdot f_2 - d_2 \cdot f_1 + \phi(d_1,d_2)) + [d_1,d_2]_{\mathcal{D}}
\end{align*}
for all $l_1, l_2 \in L, \ f_1, f_2 \in \mathcal{C}, \ d_1, d_2 \in \mathcal{D}$, where $[\cdot, \cdot]_{L}$ and $[\cdot, \cdot]_{\mathcal{D}}$ denote the bracket operations on $L$ and $\mathcal{D}$ respectively. If we now put $\mathfrak{h} := \sum_{0 \neq \alpha \in \Delta} [L_{\alpha}^0, L_{-\alpha}^0]$, then $E$ contains the ad-diagonalizable subalgebra $H = \mathfrak{h} \oplus \mathcal{C}_0 \oplus \mathcal{D}_0$ and also has an invariant symmetric non-degenerate bilinear form 
\begin{align*}
\langle l_1 + f_1 + d_1 | l_2 + f_2 + d_2 \rangle := \langle l_1,l_2 \rangle + f_1(d_2) + f_2(d_1). 
\end{align*}

\bthm \cite[Theorem 6]{Ne} \label{T3.1}
\
\begin{enumerate}
\item The Lie algebra $E(L, \mathcal{D}, \phi)$ constructed in (\ref{(3.2)}) is an EALA of nullity $n+1$, relative to the form $\langle \cdot|\cdot \rangle$ and the ad-diagonalizable subalgebra $H$. In this case, the core and the centerless core of $E$ are given by $L \oplus \mathcal{C}$ and $L$ respectively.

\item Conversely let $E$ be an EALA of nullity $n+1$ and $L = E_c/Z(E_c)$ be its centerless core. Then $L$ is a centerless Lie torus of nullity $n+1$ and there exists a unique subalgebra $\mathcal{D}$ of $\mathrm{SCDer}(L)$ which induces the $\Lambda$-grading of $L$ and a $2$-cocycle $\phi : \mathcal{D} \times \mathcal{D} \longrightarrow \mathcal{C}$ satisfying (\ref{(3.1)}) such that $E \cong E(L, \mathcal{D}, \phi)$.
\end{enumerate} 
\ethm 

\brmk
If $E$ is an EALA of nullity $0$, then $E \cong E(L, \mathcal{D}, \phi)$, where $L$ is a finite-dimensional simple Lie algebra, $\mathcal{D} = \mathcal{C} = (0)$ and $\phi = 0$. In this case, the Cartan subalgebra of $L$ plays the role of $H$ and the invariant, non-degenerate, symmetric bilinear form on $L$ is given by $(\cdot|\cdot)$. 
\ermk

\subsection{Realization of fgc EALAs} \label{SS3.3}
Let $E := E_{n+1} (n \geqslant 0)$ be an EALA of nullity $n+1$, whose centerless core $E_{cc}$ satisfies the \textit{fgc condition}.
By \cite[Proposition 3.25]{ABFP} and \cite[Theorem 3.1]{ABFP}, $E_{cc}$ is isomorphic to a Lie torus $L(\mathfrak{g}, \underline{\sigma})$ (see \S \ref{SS2.3}). We shall refer to these EALAs as \textit{ twisted fgc EALAs}. In particular, if $E_{cc} \cong L(\mathfrak{g})$ (i.e. all the automorphisms $\sigma_0, \ldots, \sigma_n$ used in our construction of $L(g, \underline{\sigma})$ in \S \ref{SS2.3} are taken to be the identity), then $E$ is called an (untwisted) \textit{fgc EALA}. In this paper, we shall only consider (untwisted) fgc EALAs. 

Let $L = L(\mathfrak{g})$ and $\Lambda = \mathbb{Z}^{n+1}$. Then $L$ is $\Lambda$-graded with Cent($L$) $\cong A_{n+1}$ \cite{BN,GP}. This isomorphism allows us to consider $L$ as an algebra over $A_{n+1}$. Using the invariant, non-degenerate form $(\cdot | \cdot)$ on $\mathfrak{g}$, we can now define an invariant, symmetric and non-degenerate bilinear form $\langle \cdot, \cdot \rangle$ on $L$ via 
\begin{align*}
\langle x \otimes t^{\underline{r}}, y \otimes t^{\underline{s}} \rangle := (x | y) \delta_{\underline{r} + \underline{s}, 0} \ \forall \ \underline{r}, \underline{s} \in \mathbb{Z}^{n+1}. 
\end{align*}
Moreover it follows from \cite[Corollary 7.4]{NPPS} that such a  form  on $L$ is also unique up to scalars.

The $\Lambda$-grading on $A_{n+1}$ induces a \textit{degree derivation} $\partial_{\theta}$ of $A_{n+1}$, given by $\partial_{\theta}(t^{\underline{k}}) = \theta(\underline{k})t^{\underline{k}}$ for any $\theta \in \mathrm{Hom}_{\mathbb{Z}}(\Lambda, \mathbb{C})$. Putting $D = \{\partial_{\theta} \ | \ \theta \in \mathrm{Hom}_{\mathbb{Z}}(\Lambda,\mathbb{C}) \}$, we then have a vector space isomorphism from $\mathrm{Hom}_{\mathbb{Z}}(\Lambda, \mathbb{C})$ to $D$ via the map $\theta \mapsto \partial_{\theta}$ and thus $D \cong \mathbb{C}^{n+1}$. This gives $\mathrm{Der}(A_{n+1}) = \mathrm{span} \{t^{\underline{r}}d_i \ | \ \underline{r} \in \mathbb{Z}^{n+1}, \ 0 \leqslant i \leqslant n \}$, where $d_i = t_i (\partial/\partial t_i)$. It is well-known that $\mathrm{Der}(A_{n+1})$ forms a $\mathbb{Z}^{n+1}$-graded Lie algebra under the bracket operation $[t^{\underline{r}}d_i,t^{\underline{s}}d_j]:=s_it^{\underline{r} + \underline{s}}d_j-r_j t^{\underline{r} + \underline{s}}d_i$.

In order to explicitly construct fgc EALAs, we need to define the subalgebra of skew-symmetric derivations with respect to $\langle \cdot,\cdot \rangle$. To this end, put
\begin{align*}
\mathcal{S}_{n+1} := \{\delta \in \mathrm{Der}(A_{n+1}) \ | \ \langle \delta(x \otimes t^{\underline{r}}) , y \otimes t^{\underline{s}} \rangle + \langle x \otimes t^{\underline{r}} , \delta(y \otimes t^{\underline{s}}) \rangle = 0 \ \forall \ x, y \in \mathfrak{g}, \ \underline{r}, \underline{s} \in \mathbb{Z}^{n+1} \} \\
= \{\sum_{i=0}^{n}u_it^{\underline{r}}d_i \in \mathrm{Der}(A_{n+1}) \ | \ \sum_{i=0}^{n}u_ir_i = 0 \} = \bigoplus_{\underline{r} \in \mathbb{Z}^{n+1}} \mathbb{C}t^{\underline{r}}\{\sum_{i=0}^{n}u_id_i \ | \ \sum_{i=0}^{n}u_ir_i = 0, \ \underline{u} \in \mathbb{C}^{n+1} \},
\end{align*}
which is the subalgebra of skew-symmetric derivations of $A_{n+1}$. It can be shown that $\mathrm{Der}(A_{n+1})$ can be identified with $\mathrm{CDer}(L)$ and under this identification, the restriction to $\mathcal{S}_{n+1}$ again gets mapped onto $\mathrm{SCDer}(L)$ \cite{AP}. We now consider any graded subalgebra $\mathcal{D} \subseteq \mathrm{SCDer}(L) \cong \mathcal{S}_{n+1}$ such that the evaluation map ev: $\Lambda \longrightarrow \mathcal{D}_{\underline{0}}^*$ is injective (see \S \ref{SS3.2}), where $\mathcal{D}_0 = \mathrm{span} \{d_i \ | \ 0 \leqslant i \leqslant n \} \cong D$. This induces the central 2-cocycle $\sigma_{\mathcal{D}}$ for $L$ with values in the graded dual $\mathcal{C}$, from which we obtain a (non-zero) central extension $L \oplus \mathcal{C}$ of $L$. Then it follows from Theorem \ref{T3.1} that $E(L, \mathcal{D}, \phi)$ is an EALA with centerless core $L$ for any affine 2-cocycle $\phi : \mathcal{D} \times \mathcal{D} \longrightarrow \mathcal{C}$ and conversely, every EALA $E$ with centerless core $L(\mathfrak{g})$ is isomorphic to $E(L(\mathfrak{g}), \mathcal{D}, \phi)$ for suitable choices of $\mathcal{D}$ and $\phi$. In this case, $H = \mathfrak{h} \oplus \mathcal{C}_{\underline{0}} \oplus \mathcal{D}_{\underline{0}}$ is an ad-diagonalizable subalgebra of $E$. Thus we have the following result.

\bcrlre \label{C3.3}
Let $E_{n+1}$ be an fgc EALA of nullity $n+1$, where $n \geqslant 0$. Then $E_{n+1} \cong L(\mathfrak{g}) \oplus \mathcal{C}_{n+1} \oplus \mathcal{D}_{n+1}$ for some (non-zero) central extension $L(\mathfrak{g}) \oplus \mathcal{C}_{n+1}$ of $L(\mathfrak{g})$ and a graded subalgebra $\mathcal{D}_{n+1}$ of $\mathcal{S}_{n+1}$ containing $(\mathcal{D}_{n+1})_{\underline{0}} = \mathrm{span} \{d_i \ | \ 0 \leqslant i \leqslant n \}$, with $\mathcal{C}_{n+1}$ being the graded dual of $\mathcal{D}_{n+1}$ and $H = \mathfrak{h} \oplus (\mathcal{C}_{n+1})_{\underline{0}} \oplus (\mathcal{D}_{n+1})_{\underline{0}}$ being an ad-diagonalizable subalgebra of $E_{n+1}$.  
\ecrlre

\brmk
\
\begin{enumerate}
\item Note that $(\mathcal{D}_{n+1})_{\underline{0}}$ is an \textit{abelian ad-diagonalizable subalgebra} of $E_{n+1}$. 
\item For $n = 0$, it is easy to see that $\mathcal{D}_1 = (\mathcal{D}_1)_0 = \mathbb{C}d = \mathcal{S}_1$, where $d = t(d/dt)$ and $\mathcal{C}_1 = \mathbb{C}K$. But this is far from true for $n \geqslant 1$, as there are plenty of choices for graded subalgebras of $\mathcal{S}_{n+1}$ in this case, which in turn give rise to different examples of fgc EALAs.
\end{enumerate}
\ermk

\subsection{Cores of fgc EALAs as quotients of the universal central extension of $L(\mathfrak{g})$} \label{SS 3.4}
Let us consider the module of differentials ($\Omega_{A_{n+1}},d$) of $A_{n+1}$, which is the
free $A_{n+1}$-module with basis $\{K_0, \ldots, K_n \}$ along with the differential map $d : A_{n+1} \longrightarrow \Omega_{A_{n+1}}$. The image of this
map is spanned by
$d(t^{\underline{k}}) = \sum_{i =0}^{n}k_i t^{\underline{k}}K_i$ for
$\underline{k} \in \ZZ^{n+1}$, where $K_i = t_i^{-1}dt_i \
\forall \ 0 \leqslant i \leqslant n$. More precisely,
\begin{align*}
\Omega_{A_{n+1}} = \text{span} \{t^{\underline{k}}K_i \ | \ 0 \leqslant i \leqslant n, \ \underline{k} \in \ZZ^{n+1} \},\ 
dA_{n+1} =  \text{span} \big \{\sum_{i=0}^{n}{k_i t^{\underline{k}}K_i} \ | \ \underline{k} \in \ZZ^{n+1} \big \}.
\end{align*}
If we now consider the quotient space $\mathcal{Z}_{n+1} = \Omega_{A_{n+1}} / dA_{n+1}$, then we know that
\begin{align}\label{Universal}
\overline{L}(\mathfrak{g})=L(\mathfrak{g}) \oplus \mathcal{Z}_{n+1}
\end{align} 
is the \textit{universal central extension} of $L(\mathfrak{g})$ \cite{Ka,EMY}. By
abuse of notation, we shall denote the image of $t^{\underline{k}}K_i$ in
$\mathcal{Z}_{n+1}$ again by itself and define the bracket operation
on $\overline{L}(\mathfrak{g})$ as follows:
\begin{enumerate}
	\item $[x \otimes t^{\underline{k}}, y \otimes t^{\underline{l}}] = [x,y] \otimes t^{\underline{k} + \underline{l}} + (x|y) \sum_{i=0}^{n} k_it^{\underline{k} + \underline{l}}K_i$;
	\item $\mathcal{Z}_{n+1}$ is central in $\overline{L}(\mathfrak{g})$.  
\end{enumerate}

\brmk
We normalize $(\cdot|\cdot)$ such that $(\theta|\theta)=2$, where $\theta$ is the highest root of $\mathfrak{g}$.
\ermk

\blmma \label{L3.6}
Let $E_{n+1}$ be an fgc EALA of nullity $n+1$, where $n \geqslant 0$. Then 
\begin{enumerate}
\item $(E_{n+1})_c \cong L(\mathfrak{g}) \oplus (\Omega_{A_{n+1}}/dA_{n+1}^{\prime})$, where
$dA_{n+1}^{\prime}$ is a subspace of $\Omega_{A_{n+1}}$ containing $dA_{n+1}$. 
\item $Z(E_{n+1}) = \mathrm{span} \{K_0, \ldots, K_n \} = (\Omega_{A_{n+1}}/dA_{n+1}^{\prime})_{\underline{0}}$. 
\end{enumerate}
\elmma

\begin{proof} 
If $n = 0$, then it follows from \cite{ABGP} that $(E_1)_c \cong g \otimes \mathbb{C}[t, t^{-1}] \oplus \mathbb{C}K$, where $ \mathcal{Z}_{1} = \Omega_{A_1}/dA_1 = \mathbb{C}K$ and so we are done. Assume that $n \geqslant 1$. By Corollary \ref{C3.3}, we have $(E_{n+1})_c = L(\mathfrak{g}) \oplus \mathcal{C}_{n+1}$ for some non-zero central extension $L(\mathfrak{g}) \oplus \mathcal{C}_{n+1}$ of $L(\mathfrak{g})$. \\ 
(1) We have the two short exact sequences: 
\begin{align*} 
0 \longrightarrow \mathcal{C}_{n+1} \longrightarrow (E_{n+1})_c {\overset{\pi} \longrightarrow} L(\mathfrak{g}) \longrightarrow 0, \\
0 \longrightarrow \mathcal{Z}_{n+1} \longrightarrow \overline{L}(\mathfrak{g}) {\overset{\pi^{\prime}} \longrightarrow} L(\mathfrak{g}) \longrightarrow 0. 
\end{align*}
Now since $\overline{L}(\mathfrak{g})$ is the universal central extension of $L(\mathfrak{g})$, there exists a unique homomorphism (of Lie algebras) $F : \overline{L}(\mathfrak{g}) \longrightarrow (E_{n+1})_c$ such that $\pi \circ F = \pi^{\prime}$. We claim that $F$ is onto. \\
From (P2) of \S \ref{SS2.2}, it follows that both $(E_{n+1})_c$ and $\overline{L}(\mathfrak{g})$ are perfect Lie algebras. In particular, we have $(E_{n+1})_c = F(\overline{L}(\mathfrak{g})) + \mathrm{Ker}(\pi)$. But as $\mathrm{Ker}(\pi) \subseteq Z((E_{n+1})_c)$, this implies that $(E_{n+1})_c = [(E_{n+1})_c, (E_{n+1})_c] = [F(\overline{L}(\mathfrak{g})),  F(\overline{L}(\mathfrak{g}))] = F([\overline{L}(\mathfrak{g}), \overline{L}(\mathfrak{g})]) = F(\overline{L}(\mathfrak{g}))$. Hence the claim. \\
The above claim gives $\mathrm{Ker}(F) \cap L(\mathfrak{g}) = (0)$ and thus $(E_{n+1}) _c\cong L(\mathfrak{g}) \oplus (\mathcal{Z}_{n+1}/\mathcal{Z}_{n+1}^{\prime})$ for some subspace $\mathcal{Z}_{n+1}^{\prime}$ of $\mathcal{Z}_{n+1}$, which thereby proves the assertion by (\ref{Universal}). \\
(2) By Corollary \ref{C3.3}, $E_{n+1} \cong (E_{n+1})_c \oplus \mathcal{D}_{n+1}$ for some graded subalgebra $\mathcal{D}_{n+1}$ of $\mathcal{S}_{n+1}$ containing $(\mathcal{D}_{n+1})_{\underline{0}}$. Now since $\mathcal{D}_{n+1}$ acts on $L(\mathfrak{g})$ by derivations, it is evident that $\mathcal{D}_{n+1} \cap Z(E_{n+1}) = (0)$ and $L(\mathfrak{g}) \cap Z(E_{n+1}) = (0)$. Moreover $f_{\underline{p}} \in (\mathcal{C}_{n+1})_{\underline{p}} \cap Z(E_{n+1})$ if and only if $(d(\underline{r})\cdot f_{\underline{p}})(d(\underline{s})) = -f_{\underline{p}}([d(\underline{r}),d(\underline{s})]) = 0 \ \forall \ d(\underline{r}) \in (\mathcal{D}_{n+1})_{\underline{r}}, \ d(\underline{s}) \in (\mathcal{D}_{n+1})_{\underline{s}}$, where $\underline{r}, \underline{s} \in \mathbb{Z}^{n+1}$ are arbitrary. This implies that $f_{\underline{p}} \in (\mathcal{C}_{n+1})_{\underline{p}} \cap Z(E_{n+1})$ if and only if $\underline{p} = \underline{0}$, which finally gives the desired result.     
\end{proof}

\brmk
As in the case of $\mathcal{Z}_{n+1}$, we shall again denote the image of $t^{\underline{k}}K_i$ in $\mathcal{C}_{n+1}$ by itself. 
\ermk
	        
\section{Examples of fgc EALAs and related Lie algebras} \label{S4}
In this section, we give concrete examples of Lie algebras which are fgc EALAs or intimately connected to such EALAs and whose representations have been previously studied in the literature.
\subsection{Full toroidal Lie algebra} \label{SS4.1}
The \textit{full toroidal Lie algebra} $\tau_{n+1}^{\mathcal{F}}$, which is a higher-dimensional analogue of the affine-Virasoro algebra, is defined to be the semi-direct product of $\overline{L}(\mathfrak{g})$ (see (\ref{Universal})) with $\mathrm{Der}(A_{n+1}) = \mathrm{span} \{t^{\underline{r}}d_i : \underline{r} \in \mathbb{Z}^{n+1}, \ 0 \leqslant i \leqslant n \}$. It is well-known that $\mathrm{Der}(A_{n+1})$ admits the following non-trivial
2-cocycles $\phi_1$ and $\phi_2$ with values in $\mathcal{Z}_{n+1}$:
\begin{equation} \label{2-Cocycle}
	\phi_1(t^{\underline{r}}d_i, t^{\underline{s}}d_j) = -s_ir_j
	\sum_{p=0}^{n}{r_p t^{\underline{r} + \underline{s}}K_p}, \qquad
	\phi_2(t^{\underline{r}}d_i, t^{\underline{s}}d_j) = r_is_j \sum_{p=0}^{n}{r_p t^{\underline{r} + \underline{s}}K_p}          
\end{equation}
(see \cite{YB} for more details).
Let $\phi$ be an arbitrary linear combination of $\phi_1$ and $\phi_2$. Then we can define the \textit{full toroidal Lie algebra} in $n+1$ variables (relative to $\mathfrak{g}$ and $\phi$) by setting
\begin{align}\label{Full Toroidal}
	\tau_{n+1}^{\mathcal{F}} := L(\mathfrak{g}) \oplus \mathcal{Z}_{n+1} \oplus \mathrm{Der}(A_{n+1}).
\end{align}
with the following bracket operations besides the relations (1) and
(2) mentioned in (\ref{Universal}): 
\begin{enumerate}
	\item $[t^{\underline{r}}d_i,t^{\underline{s}}K_j]=s_it^{\underline{r} + \underline{s}}K_j + \delta_{ij} \sum_{p=0}^{n}r_p t^{\underline{r} + \underline{s}}K_p$,
	\item $[t^{\underline{r}}d_i,t^{\underline{s}}d_j]=s_it^{\underline{r} + \underline{s}}d_j-r_j t^{\underline{r} + \underline{s}}d_i + \phi(t^{\underline{r}}d_i,t^{\underline{s}}d_j)$,
	\item $[t^{\underline{r}}d_i,x \otimes t^{\underline{s}}]=s_ix \otimes t^{\underline{r} + \underline{s}} \ \forall \ x \in \mathfrak{g},\ \underline{r},\ \underline{s} \in \mathbb{Z}^{n+1}, \ 0 \leqslant i,j \leqslant n$.   
\end{enumerate}

\brmk
\
\begin{enumerate}
\item The irreducible quasi-finite modules for the affine-Virasoro algebra were classified in \cite{CLW,LPX}. 
\item The classification of irreducible integrable representations of $\tau_{n+1}^{\mathcal{F}}$ was provided in \cite{EJ}, whereas a large class of irreducible quasi-finite modules over $\tau_{n+1}^{\mathcal{F}}$ were explicitly constructed in \cite{B} by means of vertex algebras.
\end{enumerate}     
\ermk
\subsection{Toroidal Lie algebra} \label{SS4.2}
The \textit{toroidal Lie algebra} $\tau_{n+1}^{tor}$ is formed by simply adding the space of degree derivations $D = \sum_{i=0}^{n} \mathbb{C}d_i$ to $\overline{L}(\mathfrak{g})$, i.e.
\begin{align*}
\tau_{n+1}^{tor} := L(\mathfrak{g}) \oplus \mathcal{Z}_{n+1} \oplus D. 
\end{align*}
The bracket operations on $\tau_{n+1}^{tor}$ are given by restricting the bilinear operations endowed on $\tau_{n+1}^{\mathcal{F}}$.

\brmk
Both $\tau_{n+1}^{tor}$ ($n \geqslant 1$) and $\tau_{n+1}^{\mathcal{F}}$ fall short of an EALA as they do not possess an invariant symmetric non-degenerate bilinear form. Also note that $Z(\tau_{n+1}^F) = Z(\tau_{n+1}^{tor}) = \mathrm{span} \{K_0, \ldots, K_n \}$.  
\ermk

\subsection{Toroidal extended affine Lie algebra} \label{SS4.3}
The \textit{toroidal extended affine Lie algebra} $\tau_{n+1}^{\mathcal{S}}$ of nullity $n+1$ is given by the semidirect sum of $\overline{L}(\mathfrak{g})$ and $\mathcal{S}_{n+1}$ (see \S \ref{SS3.3}), i.e.
\begin{align*}
\tau_{n+1}^{\mathcal{S}} := L(\mathfrak{g}) \oplus \mathcal{Z}_{n+1} \oplus \mathcal{S}_{n+1}, 
\end{align*}
with its bracket operations being induced from that of $\tau_{n+1}^{\mathcal{F}}$. 

\subsection{Minimal extended affine Lie algebra} \label{SS4.4}
The \textit{minimal extended affine Lie algebra} $\tau_{n+1}^{M}$ of nullity $n+1$ is formed by adjoining the space of degree derivations $D$ (which is the smallest possible subalgebra of $\mathcal{S}_{n+1}$ allowed in the construction of untwisted fgc EALAs) to the minimal (non-zero) central extension of $L(\mathfrak{g})$. More precisely,
\begin{align*}
\tau_{n+1}^{M} := L(\mathfrak{g}) \oplus \sum_{i=0}^{n} \mathbb{C}K_i \oplus D, 
\end{align*} 
which forms a Lie algebra under the following bracket operations.
\begin{enumerate}
	\item $[x \otimes t^{\underline{r}}, y \otimes t^{\underline{s}}] = [x,y] \otimes t^{\underline{r} + \underline{s}} \ + \ \delta_{\underline{r} + \underline{s}, \underline{0}}(x|y)\displaystyle{\sum_{i=0}^{n}} r_iK_i$ ;
	\item $K_i$'s are central in $\tau_{n+1}^{M}$ ;
	\item $[d_i,x \otimes t^{\underline{r}}] = r_ix \otimes t^{\underline{r}} \ \forall \ x, y \in \mathfrak{g}, \ \underline{r}, \underline{s} \in \mathbb{Z}^{n+1}, \ 0 \leqslant i \leqslant n$. 
\end{enumerate}

\subsection{Hamiltonian extended affine Lie algebra} \label{SS4.5}
The \textit{Hamiltonian extended affine Lie algebra} $\tau_{n+1}^{\mathcal{H}}$ of nullity $n+1$ (where $n = 2N + 1 \geqslant 1$), which was first introduced in \cite{E4}, is obtained by adjoining $L(\mathfrak{g})$ with the Hamiltonian Lie algebra and its graded dual. We now explicitly describe the construction of this EALA following \cite{E4,T}. For each $\underline{r} = (r_0, \ldots, r_N, r_{N+1}, \ldots, r_n) \in \mathbb{Z}^{n+1}$, set $h_{\underline{r}} = \sum_{i=0}^{N} \big(r_{N+i}t^{\underline{r}}d_i - r_it^{\underline{r}}d_{N+i} \big)$. Then the corresponding Hamiltonian Lie algebra is given by 
\begin{align*}
\mathcal{H}_{n+1} := \mathrm{span} \{h_{\underline{r}},  d_i \ | \ \underline{0} \neq \underline{r} \in \mathbb{Z}^{n+1}, \ 0 \leqslant i \leqslant n \}, 
\end{align*}
where the bracket operations are induced from $\mathcal{S}_{n+1}$. Note that $\mathcal{H}_{n+1}$ is a $\mathbb{Z}^{n+1}$-graded Lie algebra with $\mathcal{H}_{n+1} = D \ltimes \big (\bigoplus_{\underline{0} \neq \underline{r} \in \mathbb{Z}^{n+1}} (\mathcal{H}_{n+1})_{\underline{r}} \big )$ and $(\mathcal{H}_{n+1})_{\underline{0}} = D$ (see \cite{T} for more details). 

In order to determine the graded dual of $\mathcal{H}_{n+1}$, let us first put
\begin{align*}
\mathcal{K}_{n+1} := \{\sum_{i=0}^{n}u_it^{\underline{r}}K_i \in \mathcal{Z}_{n+1} \ | \ \sum_{i=0}^{n}(u_ir_{N+i} - u_{N+i}r_i) = 0 \}.
\end{align*}
It is easy to check that $[\mathcal{H}_{n+1}, \mathcal{K}_{n+1}] \subseteq \mathcal{K}_{n+1}$. Then we can define the Hamiltonian EALA as 
\begin{align*}
\tau_{n+1}^{\mathcal{H}} := L(\mathfrak{g}) \oplus \mathcal{Z}_{n+1}/\mathcal{K}_{n+1} \oplus \mathcal{H}_{n+1}. 
\end{align*} 
The Lie algebra structure on $\tau_{n+1}^{\mathcal{H}}$ is obtained by restricting the bracket operations on $\tau_{n+1}^{\mathcal{F}}$. 

\brmk
Note that $\mathcal{H}_{n+1}$ is a subalgebra of $\mathcal{S}_{n+1}$. In particular, it can be easily verified that $\mathcal{H}_2 \cong \mathcal{S}_2$, due to which we finally obtain $\tau_2^{\mathcal{H}} \cong \tau_2^{\mathcal{S}}$.
\ermk
 
\section{Roots, co-roots and Weyl group of fgc EALAs} \label{S5}

\subsection{Roots and co-roots}
Let $E_{n+1} = L(\mathfrak{g}) \oplus \mathcal{C}_{n+1} \oplus \mathcal{D}_{n+1}$ be an fgc EALA of nullity $n+1$ with an ad-diagonalizable subalgebra $H = \mathfrak{h} \oplus \displaystyle{\sum_{i =0}^{n}\CC K_i \oplus \sum_{i=0}^{n}} \CC d_i$. To describe the roots of $E_{n+1}$, let us define $\delta_i \in H^*$ by setting
\begin{align*}
\delta_i (\mathfrak{h}) =0,\,\,\,\  \delta_i (K_j) = 0\,\,\,\ \text{and}\,\,\,\ \delta_i (d_j) = \delta_{i j}\,\,\,\ \forall \ 0 \leqslant i,j \leqslant n.
\end{align*}
Put $\delta_{\underline{\gamma}} = \displaystyle{\sum_{i =0}^{n}{\gamma_i \delta_i}}$ for $\underline{\gamma} \in \mathbb{C}^{n+1}$. Then the roots of $E_{n+1}$ with respect to $H$ are given by $R = R^{\times} \sqcup R^0$, where $R^{\times} = \{ \alpha + \delta_{\underline{k}}\ | \ \alpha \in \Delta_{\mathfrak{g}}^{\times}, \ \underline{k} \in \ZZ^{n+1} \}$ and $R^0 = \{\delta_{\underline{k}} \ | \ \underline{k} \in \mathbb{Z}^{n+1} \}$. Also for each $\beta = \alpha + \delta_{\underline{k}} \in R^{\times}$, define the corresponding co-root $\beta^{\vee} := \alpha^{\vee} + \frac{2}{(\alpha|\alpha)} \displaystyle{\sum_{i =0}^{n}{k_i K_i}}$ where $\alpha^{\vee} \in \mathfrak{h}$ is the co-root of $\alpha \in \Delta_{\mathfrak{g}}^{\times}$.       

\subsection{The Weyl group}
For each $\beta \in R^{\times}$, define
the reflection operator $r_{\beta}$ on $H^*$ by setting
\begin{align*}
r_{\beta}(\lambda) := \lambda - \lambda(\beta^{\vee}) \beta \ \forall \ \lambda \in H^*.
\end{align*}
Then the Weyl group of $E_{n+1}$, which we shall denote by $\mathcal{W}$, is the group
generated by all such reflections $r_{\beta}$ with $\beta \in R^{\times}$. Note that the Weyl group of $E_1$ is the affine Weyl group. 

\section{Automorphism twist} \label{S6} 
Let $E_{n+1}$ be an fgc EALA of nullity $n+1$ and $G=GL(n+1,\mathbb{Z})$. Then $G$ acts on $\mathbb{Z}^{n+1}$ by means of matrix multiplication. Let us now fix any $A = (a_{pj})_{0 \leqslant p,j \leqslant n}$ in $G$ and define 
\[
T_{A}(x \otimes t^{\underline{m}}) = x \otimes t^{\underline{m}A^t},
\qquad
T_{A}(t^{\underline{m}}K_j) = \sum_{p=0}^{n}a_{pj}
t^{\underline{m}A^t}K_p, \qquad
T_{A}(t^{\underline{m}}d_j) = \sum_{p=0}^{n}b_{jp}
t^{\underline{m}A^t}d_p,
\]
where $0 \leqslant j \leqslant n$, $B = (b_{jp}) = A^{-1}, \ A^{t}$ is the transpose of $A$ and $\underline{m}$ is a row vector in $\mathbb{Z}^{n+1}$. One can check that $T_{A}$ is an automorphism of $\tau_{n+1}^{\mathcal{F}}$, which leaves  $\tau_{n+1}^{\mathcal{S}}, \ \tau_{n+1}^{M}$ and $\tau_{n+1}^{tor}$ invariant. In the general set-up, it follows from Corollary \ref{C3.3} and Lemma \ref{L3.6} that $T_B$ takes $E_{n+1}$ to a Lie algebra $(E_{n+1})^A$ satisfying $E_{n+1} \cong (E_{n+1})^A$ and such a phenomenon is also true for the Lie algebra $L(\mathfrak{g}) \oplus \mathcal{Z}_{n+1} \oplus \sum_{i=0}^{n} \mathbb{C}d_i$. For notational convenience, we shall use this identification to (loosely) say that $G$ acts on $E_{n+1}$ via automorphisms. In the current paper, we shall use this notion without any further comments and simply refer to it as \textit{up to a twist of an automorphism}.

\brmk
The above discussion reveals that whenever we twist an $E_{n+1}$-module $V$ by $A \in G$, the resulting module will be a module over a Lie algebra isomorphic to $E_{n+1}$. Nonetheless, to avoid notational complexity, we shall again consider the twisted module over the original Lie algebra. 
\ermk

\section{Some preliminary definitions and results} \label{S7}
In this section, we introduce the notion of quasi-finite modules of an fgc EALA $E := E_{n+1} = L(\mathfrak{g}) \oplus \mathcal{C} \oplus \mathcal{D}$ of nullity $n+1 (\geqslant 1)$ and prove some general results related to these modules, which will be  utilized in the subsequent sections.   

\bdfn
$V$ is said to be a quasi-finite module over $E$ if it
satisfies:
\begin{enumerate}
	\item $V = \bigoplus_{\underline{r} \in \mathbb{Z}^{n+1}} V_{\underline{r}}, \ \text{where} \ V_{\underline{r}} = \{v \in V \ | \ d_iv = r_iv, \ 0 \leqslant i \leqslant n \}$;
	\item $\dim V_{\underline{r}} < \infty \ \forall \ \underline{r} \in \mathbb{Z}^{n+1}$.
\end{enumerate}  
\edfn
\noindent The collection $P_{D}(V) := \{\mu \in D^* \ | \ V_{\mu} \neq (0) \}$, where $V_{\mu} := \{v \in V \ | \ d_iv = \mu(d_i)v, \ 0 \leqslant i \leqslant n \}$ and $D= \mathrm{span} \{d_0, \ldots, d_n \}$, is known as the set of all $D$-weights of $V$. Also set $\tau_{n+1} := L(\mathfrak{g}) \oplus \mathcal{C} \oplus D$.

\bdfn
A quasi-finite module $V$ over $E$ is said to be uniformly bounded if there exists $N_0 \in \mathbb{N}$ such that $\mathrm{dim} V_{\underline{r}} \leqslant N_0 \ \forall \ \underline{r} \in \mathbb{Z}^{n+1}$.
\edfn

\brmk\label{R7.3}
\
\begin{enumerate}
\item If $V$ is an irreducible (and hence indecomposable) quasi-finite module over $E_{n+1}$, then it is easy to see that there exists $\lambda \in D^*$ such that $P_{D}(V) \subseteq \{\lambda + \underline{r} \ | \ \underline{r} \in \mathbb{Z}^{n+1} \}$. 
\item By Lemma \ref{L3.6}, each $K_i \in \mathcal{C}_{\underline{0}}$ acts by a fixed scalar on $V$, say $c_i$. If $c_i = 0 \ \forall \ 0 \leqslant i \leqslant n$, then we say that $V$ is a \textit{level zero module}, otherwise we say that $V$ has \textit{non-zero level}. 
\item For a quasi-finite module $V$ over $E$,  $V_{\underline{r}}$ is a finite-dimensional $\mathfrak{g}$-module and thus $V_{\underline{r}}$ has a weight space decomposition with respect to $\mathfrak{h}$ for each $\underline{r} \in \mathbb{Z}^{n+1}$. Consequently
\begin{align*}
V = \bigoplus_{\eta \in \mathfrak{h}^*, \ \underline{r} \in \mathbb{Z}^{n+1}} V_{\eta + \delta_{\underline{r}}}, \text{where} \,\ V_{\eta + \delta_{\underline{r}}} = V_{\eta} \cap V_{\underline{r}}, \,\ V_{\eta} = \{v \in V \ | \ hv=\eta(h)v \ \forall \ h \in \mathfrak{h} \}.
\end{align*} 
This induces an $H$-weight space decomposition of $V$ and the set of all $H$-weights of $V$ is given by $P_{H}(V) := \{\nu \in H^* \ | \ V_{\nu} \neq (0) \}$, where $V_{\nu} := \{v \in V \ | \ hv = \nu(h)v \ \forall \ h \in H \}$ is finite-dimensional for each $\nu \in H^*$.   
\item For an $E$-module $V$, $\{v \in V \ | \ (E_c)v = 0 \}$ is an $E$-submodule of $V$, as $E_c$ is an ideal of $E$.
\end{enumerate}  
\ermk 

\blmma \label{L7.4}
If $V$ is an irreducible quasi-finite module over $E$ of non-zero level, then up to a twist of an automorphism, we can assume that $c_0 \neq 0$ and $c_1 = \ldots = c_n = 0$. 
\elmma

\begin{proof}
Pick $0 \neq v \in V$. By Zorn's Lemma, $U(\tau_{n+1})v$ admits a non-zero irreducible $\tau_{n+1}$-quotient, say $V^{\prime}$, with finite-dimensional $H$-weight spaces. Then we can apply \cite[Theorem 4.5]{E2} along with \cite[Theorem 1.10]{E3} and Lemma \ref{L3.6} to deduce that, up to a twist of an automorphism, $K_i$ acts trivially on $V^{\prime}$ for all $1 \leqslant i \leqslant n$. This proves the lemma, as each $K_i$ acts by a fixed scalar on $V$.
\end{proof}

\bdfn
An $E$-module $V$ is called integrable if 
\begin{enumerate}
	\item $V$ is a $H$-weight module, i.e. $\displaystyle{V = \bigoplus_{\nu \in H^*} {V_{\nu}}}$, where $V_{\nu} = \{v \in V |\,\, hv = \nu(h)v \,\, \forall \,\, h \in H \}$;
	\item The $H$-weight spaces of $V$ are finite-dimensional, i.e. dim$V_{\nu} < \infty \ \forall \ \nu \in H^*$;
	\item For each $x_{\alpha} \otimes t^{\underline{k}} \in \mathfrak{g}_{\alpha}\otimes \mathbb{C}t^{\underline{k}}$ $(\alpha \neq 0, \ \underline{k} \in \mathbb{Z}^{n+1})$ and every $v \in V$, there exists some $m=m(\alpha,\underline{k},v) \in \mathbb{N}$ such that $(x_{\alpha} \otimes t^{\underline{k}})^mv=0$.      
\end{enumerate} 
\edfn

\blmma \cite[Lemma 2.3]{E2} \label{L7.6}
Let $V$ be an integrable (not necessarily irreducible) $E$-module. Then
\begin{enumerate}
	\item $P_H(V)$ is invariant under the action of the Weyl group $\mathcal{W}$.
	\item $\mathrm{dim}(V_{\nu}) = \mathrm{dim}(V_{w \nu})$ $\forall\ \nu \in P_H(V)$ and $w \in \mathcal{W}$.
	\item If $\nu \in P_H(V)$  and $\beta \in R^{\times}$, then $\nu (\beta^{\vee}) \in \ZZ$.
	\item If $\nu \in P_H(V)$  and $\beta \in R^{\times}$ with $\nu(\beta^{\vee}) >0$, then $\nu - \beta \in P_H(V)$.
\end{enumerate}
\elmma

The following definition is motivated from the works of Kac--Wakimoto \cite{KW} and Rao--Futorny \cite{EF}, where they studied the so-called weakly integrable modules in the context of affine Lie superalgebras.
 
\bdfn
An $E$-module $V$ is said to be weakly integrable if 
\begin{enumerate}
\item $V$ has finite-dimensional $H$-weight spaces;
\item For each $v \in V$, there exists  $k=k(\alpha,v) \in \mathbb{N}$ such that $\mathfrak{g}_{\alpha}^kv = (0) \ \forall \ \alpha \in \Delta_{\mathfrak{g}}^{\times}$.      
\end{enumerate} 
\edfn

\blmma \label{L7.8}
An $E$-module $V$ is weakly integrable if and only if $V$ is a quasi-finite module. 
\elmma

\begin{proof}
By hypothesis and Remark \ref{R7.3}, $V$ is a weight module with respect to $H$ and so we have
\begin{align*}
V = \bigoplus_{\underline{r} \in \mathbb{Z}^{n+1}} V_{\underline{r}}, \ \text{where} \ V_{\underline{r}} = \{v \in V \ | \ d_iv = r_iv \ \forall \ 0 \leqslant i \leqslant n \}.
\end{align*}
If $V$ is weakly integrable, then $V_{\underline{r}}$ is an integrable $\mathfrak{g}$-module with finite-dimensional $\mathfrak{h}$-weight spaces. This conveys that $V_{\underline{r}}$ is finite-dimensional for all $\underline{r} \in \mathbb{Z}^{n+1}$ (see \cite[Lemma 3.5]{P2}). Conversely, if $V$ is a quasi-finite module, then $\mathfrak{g}$ acts on every finite-dimensional $D$-weight space, which implies that $V$ is weakly integrable.        
\end{proof}

\section{Quasi-finite modules over fgc EALAs of nullity 1 (affine case)} \label{S8}
In this section, we completely classify the irreducible quasi-finite (equivalently weakly integrable) modules over (untwisted) affine Kac-Moody algebras.

Let $\widehat{L}(\mathfrak{g}) = \mathfrak{g} \otimes \mathbb{C}[t,t^{-1}] \oplus \mathbb{C}K \oplus \mathbb{C}d$ be the affine Kac--Moody algebra (relative to $\mathfrak{g}$) with a Cartan subalgebra $H = \mathfrak{h} \oplus \mathbb{C}K \oplus \mathbb{C}d$ and whose bracket operations are obtained by simply considering $n=0$ for the minimal EALA $\tau_{n+1}^{M}$ (see \S \ref{SS4.4}). We now recall the notion of loop modules from \cite{CP}.

For any $k \in \mathbb{N}$, $\underline{\lambda} = (\lambda_1, \ldots, \lambda_k) \in (P_{\mathfrak{g}}^{+})^{k}$ and $\underline{a} = (a_1, \ldots, a_k) \in (\mathbb{C}^{\times})^k$, $a_i \neq a_j \ \forall \ 1 \leqslant i \neq j \leqslant k$, consider the \textit{loop module} $\mathcal{L}(\underline{\lambda}, \underline{a}, k, \beta) := V(\lambda_1) \otimes \ldots \otimes V(\lambda_k) \otimes \mathbb{C}[t,t^{-1}]$. Here $V(\lambda_i)$ denotes the finite-dimensional irreducible $\mathfrak{g}$-module with highest weight $\lambda_i \in P_{\mathfrak{g}}^{+}$ for all $1 \leqslant i \leqslant k$. Fix any $\beta \in \mathbb{C}$ and define an $\widehat{L}(\mathfrak{g})$-action on $\mathcal{L}(\underline{\lambda}, \underline{a}, k, \beta)$ as follows.
\begin{align*}
(x \otimes t^r)(v_1 \otimes \ldots \otimes v_k \otimes t^s) = \sum_{i=1}^{k}a_i^{r}v_1 \otimes \ldots \otimes (xv_i) \otimes \ldots \otimes v_k \otimes t^{r+s}, &\ &\
K(v_1 \otimes \ldots \otimes v_k \otimes t^s) = 0, \\ d(v_1 \otimes \ldots \otimes v_k \otimes t^s) = (\beta + s) (v_1 \otimes \ldots \otimes v_k \otimes t^s).
\end{align*}

\bthm \cite{C,CP} \label{T8.1}
\
\begin{enumerate}
\item $\mathcal{L}(\underline{\lambda}, \underline{a}, k, \beta)$ is a completely reducible integrable $\widehat{L}(\mathfrak{g})$-module.
\item If $V$ is a level zero irreducible integrable module over $\widehat{L}(\mathfrak{g})$, then $V$ is isomorphic to an irreducible component of $\mathcal{L}(\underline{\lambda}, \underline{a}, k, \beta)$ for some $k \in \mathbb{N}$, $\beta \in \mathbb{C}$, $\underline{\lambda} = (\lambda_1, \ldots, \lambda_k) \in (P_{\mathfrak{g}}^{+})^{k}$ and $\underline{a} = (a_1, \ldots, a_k) \in (\mathbb{C}^{\times})^k$, with $a_i \neq a_j \ \forall \  i \neq j$. 
\end{enumerate}
\ethm

\brmk
We shall 
refer to the above irreducible modules as \textit{irreducible loop modules}.
\ermk
 
\blmma \label{L8.3} 
Let $V$ be an irreducible quasi-finite module over a nullity $1$ (untwisted) EALA $E_1$. Then $V$ is either a highest weight or a lowest weight module or an irreducible loop module. 
\elmma

\begin{proof} 
By hypothesis, $E := E_1 \cong \widehat{L}(\mathfrak{g})$ for some finite-dimensional simple Lie algebra $\mathfrak{g}$ \cite{ABGP}. Now if the core $E_c := (E_{1})_c$, namely $\mathfrak{g} \otimes \mathbb{C}[t,t^{-1}] \oplus \mathbb{C}K$, acts trivially on $V$, then $V$ is clearly isomorphic to a $1$-dimensional highest weight module $V(a\delta)$ for some $a \in \mathbb{C}$, where $\delta \in H^*$ satisfies $\delta(d)=1$ and $\delta(K) = \delta(h) = 0 \ \forall \ h \in \mathfrak{h}$. So suppose $E_c$ acts non-trivially on $V$. \\
\textbf{Claim 1.} $V$ is either uniformly bounded or a highest weight module or a lowest weight module. \\
Let $V$ be neither a uniformly bounded module nor a lowest weight module. Fix any $\Lambda_0 \in P(V)$. Then we can find $m \in \mathbb{N}$ such that
\begin{align*}
\dim V_{\Lambda_0 - m} > (\dim \mathfrak{g}) \dim V_{\Lambda_0} + (\dim \mathfrak{g}) \dim V_{\Lambda_0 + 1}.  
\end{align*}
Thus there exists $0 \neq v \in V_{\Lambda_0 - m}$ such that $(x \otimes t^{m})v = (x \otimes t^{m+1})v = 0 \ \forall \ x \in \mathfrak{g}$. Now by using the commutator relations on $E_c$, we can deduce that $(E_c)_{j}v = (0) \ \forall \ j \geqslant m^2$, where $(E_c)_j = \text{span}\{x \otimes  t^{j}, \ \delta_{0,j}K \ | \ x \in \mathfrak{g} \}$. Again $V$ is a $\mathbb{Z}$-graded-irreducible module over $E_c$ with finite-dimensional graded components. Thus we can conclude from \cite[Lemma 1.6]{OM} that $V$ is a highest weight module, since $E_c$ is a finitely generated Lie algebra and $[E_{c}^{-}, ({E_c})_{\geqslant k}] = E_c$ for any $k \in \mathbb{N}$, where $E_c^{-} = \bigoplus_{j \in \mathbb{N}} \mathfrak{g} \otimes \mathbb{C}t^{-j}$ and $(E_c)_{\geqslant k} = \text{span}\{x \otimes  t^{j} \ | \ x \in \mathfrak{g}, \ j \geqslant k \}$. Hence the claim. \\
\textbf{Claim 2.} If $V$ is uniformly bounded (not necessarily irreducible) over $E$, then $V$ is integrable.\\
As $V$ is uniformly bounded, $K$ acts trivially on $V$ \cite{BL}. Also in view of Remark \ref{R7.3}, $V$ has finite-dimensional $H$-weight spaces. 
We now show that only finitely many $\mathfrak{h}$-weights can possibly occur in the weight space decomposition of $V$. Clearly for any $r \in \mathbb{Z}$, there exist $\mu_1, \ldots, \mu_{p(r)} \in P_{\mathfrak{g}}^{+}$ such that $V_{r} \cong \oplus_{i=1}^{p(r)} V(\mu_{i,r})$ as $\mathfrak{g}$-modules. For every $\alpha \in \Delta_{\mathfrak{g}}^{+}$, let $\mathfrak{sl}(\alpha):= \mathrm{span}\{x_{\alpha}, \alpha^{\vee}, y_{\alpha} \} \cong \mathfrak{sl}_2$. Then $V(\mu_{i,r}(\alpha^{\vee}))$ is an irreducible $\mathfrak{sl}({\alpha})$-module with $\mathrm{dim}\big(V(\mu_{i,r}(\alpha^{\vee})) \big ) = \mu_{i,r}(\alpha^{\vee})+1 \in \mathbb{N}$. Again since $V$ is uniformly bounded, there exists $N_0 \in \mathbb{N}$ with $0 \leqslant \mu_{i,r}(\alpha^{\vee}) \leqslant N_0 \ \forall \ r \in \mathbb{Z}, \ 1 \leqslant i \leqslant p(r)$ and each $\alpha \in \Delta_{\mathfrak{g}}^{+}$. This implies that there exist only finitely many $\lambda_1, \ldots, \lambda_l \in P_{\mathfrak{g}}^{+}$ such that $V$ can be expressed as a direct sum (possibly infinitely many copies) of these $V(\lambda_i)$'s as $\mathfrak{g}$-modules.\\
Next we establish that $x_{\alpha} \otimes t^k$ acts locally nilpotently on $V$ for each $\alpha \in \Delta_{\mathfrak{g}}^{\times}, \ k \in \mathbb{Z}$. \\
If not, then $x_{\alpha} \otimes t^k$ does not act locally nilpotently on $V$ for some $\alpha \in \Delta_{\mathfrak{g}}^{\times}$ and $k \in \mathbb{Z}$. Set $\gamma = \alpha + k\delta$.      
As a result, there exists $\Lambda_0^{\prime} \in P_{H}(V)$ such that $V_{\Lambda_0^{\prime} + s\gamma} \neq (0)$ for infinitely many $s \in \mathbb{N}$. But then infinitely many $\mathfrak{h}$-weights occur in the weight space decomposition of $V$, which contradicts our previous assertion. Hence the claim. The desired result now follows from Theorem \ref{T8.1}.
\end{proof}

\subsection{Classification of weakly integrable irreducible highest/lowest weight modules} \label{SS8.1}
Let us denote the simple co-roots of $\widehat{L}(\mathfrak{g})$ by $\Pi^{\vee} = \{\alpha_i^{\vee} \}_{i=0}^{l}$, where $\alpha_0^{\vee} = K - \theta^{\vee}$. Now pick any $\Lambda \in H^*$ and consider the irreducible highest weight module $V(\Lambda)$ over $\widehat{L}(\mathfrak{g})$ with highest weight $\Lambda$.\\
\textbf{Claim.} $V(\Lambda)$ is weakly integrable if and only if $\Lambda|_{\mathfrak{h}} \in P_{\mathfrak{g}}^{+}$.\\
Indeed, if $V(\Lambda)$ is weakly integrable, then $W = U(\mathfrak{g})v_{\Lambda}$ is finite-dimensional, where $V(\Lambda)_{\Lambda} = \mathbb{C}v_{\Lambda}$. By standard arguments involving $\mathfrak{sl}_2$-theory, it is now easy to conclude that $\Lambda|_{\mathfrak{h}} \in P_{\mathfrak{g}}^{+}$. Conversely, if $\Lambda|_{\mathfrak{h}} \in P_{\mathfrak{g}}^{+}$, then $W$ is a highest weight module over $\mathfrak{g}$ with highest weight $\Lambda|_{\mathfrak{h}} \in P_{\mathfrak{g}}^{+}$ and so $\mathfrak{g}_{\alpha}$ acts locally nilpotently on $v_{\Lambda}$ for each $\alpha \in \Delta_{\mathfrak{g}}^{\times}$. Finally, as $V(\Lambda) = U(\widehat{L}(\mathfrak{g}))v_{\Lambda}$ and the adjoint representation of $\widehat{L}(\mathfrak{g})$ is integrable, the claim follows immediately from \cite[Lemma 3.4]{K}.

The irreducible highest weight modules obtained in Lemma \ref{L8.3} are with respect to the triangular decomposition $\widehat{L}(\mathfrak{g}) = (\mathfrak{g} \otimes t^{-1}\mathbb{C}[t^{-1}]) \bigoplus (\mathfrak{g} \oplus \mathbb{C}K \oplus \mathbb{C}d) \bigoplus (\mathfrak{g} \otimes t\mathbb{C}[t])$. But as these modules are quasi-finite, their highest weight spaces over $\mathfrak{g} \oplus \mathbb{C}K \oplus \mathbb{C}d$ are finite-dimensional. So in the corresponding highest weight space, we can find a non-zero vector which is killed by $\mathfrak{n}^{+}$, thus implying that any such highest weight module is isomorphic to $V(\Lambda)$ for some $\Lambda \in H^*$ with $\Lambda|_{\mathfrak{h}} \in P_{\mathfrak{g}}^{+}$. Similarly, the irreducible lowest weight modules in Lemma \ref{L8.3} are all isomorphic to the graded dual $V(\Lambda)^*$ of $V(\Lambda)$, where $\Lambda|_{\mathfrak{h}} \in P_{\mathfrak{g}}^{+}$. The next result now directly follows from Lemma \ref{L7.8} and Lemma \ref{L8.3}.    

\bthm \label{T8.5}
If $V$ is an irreducible weakly integrable $\widehat{L}(\mathfrak{g})$-module, then $V$ is either isomorphic to an irreducible loop module or a highest weight module $V(\Lambda)$ or its graded dual $V(\Lambda)^*$, with $\Lambda|_{\mathfrak{h}} \in P_{\mathfrak{g}}^{+}$.  
\ethm

\brmk \label{R8.4}
The irreducible modules in Theorem \ref{T8.5} also appear in the unpublished work of Dimitrov--Grantcharov \cite{DG}, where they announced the classification of irreducible $\widehat{L}(\mathfrak{g})$-modules with finite-dimensional $H$-weight spaces. However, to the best of our knowledge, all the  published papers related to this classification problem deal only with the following cases: (i) \textit{integrable} modules \cite{C,CP}, (ii) \textit{non-zero level} modules \cite{FT,EF} and (iii) \textit{level zero} modules whose \textit{restriction to the loop algebra is reducible} \cite{E}. But there exist irreducible weakly integrable modules which do not satisfy (i), (ii) or (iii) (and obviously occur in Theorem \ref{T8.5}), as we shall see in Example \ref{E8.6}. 
\ermk

\bexmp \label{E8.6}
Define $\Lambda \in H^*$ by setting $\Lambda(\alpha_i^{\vee}) = i \ \forall \ 1 \leqslant i \leqslant l$ and $\Lambda(K) = \Lambda(d) = 0$. By construction, we then have $\Lambda(\alpha_0^{\vee}) < 0$. This implies that $V(\Lambda)$ is a \textit{level zero non-integrable (but weakly integrable) $\widehat{L}(\mathfrak{g})$-module whose restriction to $\mathfrak{g} \otimes \mathbb{C}[t,t^{-1}]$ is irreducible} (by \cite[Lemma 9.10]{K}).
\eexmp
 
\section{Quasi-finite modules over fgc EALAs of higher nullity} \label{S9}
In this section, we provide a description of the irreducible quasi-finite modules with non-trivial core action over an fgc EALA $E := E_{n+1} = L(\mathfrak{g}) \oplus \mathcal{C} \oplus \mathcal{D}$ of nullity $n+1 (> 1)$. 

\subsection{Level zero integrable modules over fgc EALAs} \label{SS9.1}

\medskip

Throughout this subsection, $V$ will denote a level zero irreducible integrable module over $E$ of nullity at least $2$ (or over the toroidal Lie algebra $\tau_{n+1}^{tor}$ for $n > 0$). We also assume that $E_c$ (respectively $\overline{L}(\mathfrak{g})$) acts non-trivially on $V$.
 
\blmma \label{L9.1} 
There exists a non-zero $v \in V$ such that $(\mathfrak{n}^{+} \otimes A_{n+1})v = (0)$.
\elmma

\begin{proof}
By hypothesis, $V$ is an integrable module over $\tau_{n+1}$ having finite-dimensional $H$-weight spaces. The desired result now follows by proceeding exactly as in \cite[Proposition 3.7]{P2}  (also see \cite[Theorem 2.4(ii)]{C} and \cite[Lemma 2.6]{E2}). Note that the irreducibility of $V$ is not required.
\end{proof}

\brmk
The above lemma conveys that $V$ can be thought of as an irreducible highest weight module with respect to the triangular decomposition $E = E(-) \oplus E(0) \oplus E(+)$, where $E(-) = \mathfrak{n}^{-} \otimes A_{n+1}, \ E(0) = \mathfrak{h} \otimes A_{n+1} \oplus \mathcal{C} \oplus \mathcal{D}, \ E(+) = \mathfrak{n}^{+} \otimes A_{n+1}$. 
\ermk

\noindent For each $j = 0, \ldots, n$, consider the subalgebras of $E$ given by  
\begin{align*}
L_j = \mathfrak{g} \otimes \mathbb{C}[{t_j},{t_j}^{-1}] \bigoplus \mathbb{C}K_j \bigoplus \mathbb{C}d_j.
\end{align*}
Let $\mathcal{W}_j$ and $\mathcal{W}_{\mathfrak{g}}$ be the Weyl groups associated to $L_j$ and $\mathfrak{g}$ respectively. Let $\theta$ denote the highest root of $\mathfrak{g}$ and  ${\theta}^{\vee}$ be the  corresponding co-root. Then by \cite[Proposition 6.5]{K}, we have
\begin{align} \label{Weyl}
\mathcal{W}_j \cong \mathcal{W}_{\mathfrak{g}} \ltimes T_j,
\end{align}
where $T_j =\{t_{\alpha_{j}} \ | \ \alpha_j \in M_j\}$ and $M_j =\gamma_j(\mathbb{Z} [\mathcal{W}_{\mathfrak{g}}\theta^{\vee}])$,
with the natural isomorphism $\gamma_j : \mathfrak{h} \longrightarrow {\mathfrak{h}}^*$.

\blmma \label{L9.3}
\
\begin{enumerate}
	\item $V(+) = \{v \in V \ | \ (\mathfrak{n}^{+} \otimes A_{n+1})v = (0) \}$ is an irreducible module over $E(0)$.
	\item There exists a unique $\overline{\lambda} \in P_{\mathfrak{g}}^{+}$ such that $P_{D}(V(+)) \subseteq \{ \overline{\lambda} + \delta_{\underline{r}} \ | \ \underline{r}  \in \mathbb{Z}^{n+1} \}$.
\end{enumerate}
\elmma

\begin{proof}
(1) Follows using the PBW theorem and the irreducibility of $V$. \\
(2) Since $[\mathfrak{h}, E(0)] = 0$, $\mathfrak{h}$ acts by scalars on $V(+)$ by (1). As a result, $\mathfrak{h}$ acts by a single linear functional on $V(+)$, which we shall denote by $\overline{\lambda}$. Finally it follows from Lemma \ref{L7.6} that $\overline{\lambda} \in P_{\mathfrak{g}}^+$.
\end{proof}

\brmk \label{R9.4}
\
\begin{enumerate}
\item If $\overline{\lambda} = 0$ in Lemma \ref{L9.3}, then for any $\alpha \in \Delta_{\mathfrak{g}}^{+}$ and $\underline{r} \in \mathbb{Z}^{n+1}$, $\text{span} \{x_{\alpha} \otimes t^{\underline{r}}, \alpha^{\vee}, y_{\alpha} \otimes t^{-\underline{r}} \} \cong \mathfrak{sl}_2$. But as $\overline{\lambda}(\alpha^{\vee}) = 0$, it follows from $\mathfrak{sl}_2$-theory for integrable modules that $(y_{\alpha} \otimes t^{-\underline{r}})v = 0$ for any fixed $0 \neq v \in V(+)$. This implies that $L(\mathfrak{g})v = 0$, which gives $(E_c)v = 0$. But then by Remark \ref{R7.3}, it is clear that $E_c$ acts trivially on $V$, due to the irreducibility of $V$. 
\item If we set $H_j = \mathfrak{h} \bigoplus \ \mathbb{C}K_j \oplus \mathbb{C}d_j$ and define $t_{\alpha_j}$ as in (\ref{Weyl}), then as $K_j$ acts trivially on $V$, it immediately follows from \cite[Equation 6.5.5]{K} that
\begin{align} \label{Translation}
t_{\alpha_j}(\mu_j) = \mu_j - \mu_j ({\alpha_j}^{\vee}) \delta_j  \ \forall \ \mu_j \in H_j^{*},\ \alpha_j \in M_j. 
\end{align}
\end{enumerate}
\ermk

The proof of the following proposition is inspired from \cite{P1}.
\bppsn \label{P9.5}
\
\begin{enumerate}
\item $\mathcal{C}$ acts trivially on $V$.
\item $V$ is an irreducible integrable module over $\mathfrak{g} \otimes A_{n+1} \oplus \mathcal{D}$ and the associated highest weight space $V(+) = \{v \in V \ | \ (\mathfrak{n}^+ \otimes A_{n+1})v = (0) \} \neq (0)$ is a uniformly bounded irreducible module over $\mathfrak{h} \otimes A_{n+1} \oplus \mathcal{D}$, i.e. $V$ is a highest weight module induced from an irreducible uniformly bounded module over $\mathfrak{h} \otimes A_{n+1} \oplus \mathcal{D}$.    
\end{enumerate}
\eppsn

\begin{proof}
(1) For each $j = 0, \ldots n$, we have $\dfrac{2\theta}{(\theta|\theta)} = \gamma_j(\theta^{\vee}) \in M_j$. Let $\theta_j = \gamma_j(\theta^{\vee})$ and $p_j = (\overline{\lambda},{\theta_j})\in \mathbb{N}$ (by Lemma \ref{L9.3} and Remark \ref{R9.4}). Now for any $k_j \in \mathbb{Z}$  with $|k_j| \geqslant p_j$ , there exists $0 \neq q_j \in \mathbb{Z}$ and $\ r_j \in \mathbb{Z}_{+}$ such that $k_j = q_jp_j + r_j,\ |r_j| \ \textless \ |p_j|$. Thus if we take $w = (\displaystyle \prod_{q_j \textless 0} t_{\theta_j}^{-q_j}) (\prod_{q_j \textgreater 0} t_{\theta_j}^{q_j}) \in \mathcal{W}$,
	then (\ref{Translation}) gives $w(\overline{\lambda} + \delta_{\underline{k}}) = \overline{\lambda} + \delta_{\underline{r}},\ |r_i| \ \textless \ |p_i| \ \forall \ i = 0, \ldots , n$. Set $P := \{ \overline{\lambda} + \delta_{\underline{r}} \ : \underline{r} \in \mathbb{Z}^{n+1}, |r_i| \ \textless \ |p_i| \}$ and $N = \oplus_{\mu \in P} (V(+))_{\mu}$.
    Again any $\tau_{n+1}$-submodule of $V$ intersecting $V(+)$ non-trivially, clearly contains a $\tau_{n+1}$-module, generated by the elements of $V(+)$. Hence by the above argument, along with Lemma \ref{L7.6}, it follows that any such submodule consists of a $\tau_{n+1}$-submodule generated by a subset of $N$.  
    Consequently we obtain a decreasing sequence of $\tau_{n+1}$-submodules of $U(\tau_{n+1})V(+)$, given by $U(\tau_{n+1})W_1 \supseteq U(\tau_{n+1})W_2 \supseteq \ldots \ldots$, where each $W_i$ is a non-zero subspace of $N$. Take $N_i^{\prime} = U(\tau_{n+1})W_i$ for each $i \in \mathbb{N}$.\\       
	\noindent \textbf{Claim.} There exists $m \in \mathbb{N}$ such that $N_m^{\prime}$ does not contain a proper $\tau_{n+1}$-submodule of $V$ intersecting $V(+)$ non-trivially. \\
	Consider the triangular decomposition $\tau_{n+1} = \tau_{n+1}^{-}  \oplus \tau_{n+1}^0 \oplus \tau_{n+1}^{+}$, where $\tau_{n+1}^{\pm} = \mathfrak{n}^{\pm} \otimes A_{n+1}$ and $\tau_{n+1}^0 = \mathfrak{h} \otimes A_{n+1} \oplus \mathcal{C} \oplus (\sum_{i=0}^{n}\mathbb{C}d_i)$. Put
	$W_i^{\prime} = U(\tau_{n+1}^0)W_i, \ \widetilde{W_i} = W_i^{\prime} \cap N$ $\forall \ i=1,2$. 
    Then we can deduce that the decreasing sequence $N_1^{\prime} \supseteq N_2^{\prime} \supseteq \ldots \ldots$ again gives rise to a decreasing sequence of (finite-dimensional) subspaces of $N$, given by $\widetilde{W_1} \supseteq \widetilde{W_2} \supseteq \ldots \ldots$, which establishes the claim. 
	
	By our claim, there exists a minimal $\tau_{n+1}$-submodule of $V$ intersecting $V(+)$ non-trivially, say $V_{min}$. Then $V_{min}$ is not necessarily irreducible over $\tau_{n+1}$. But it admits a non-zero irreducible quotient over $\tau_{n+1}$ such that $V(+)$ goes injectively to the quotient.  
	Denote this irreducible integrable $\tau_{n+1}$-module by $\mathcal{M}$. Due to Lemma \ref{L3.6}, we can now directly appeal to \cite[Proposition 4.13]{E2} to conclude that $\mathcal{C}$ acts trivially on  $\mathcal{M}$ and so $\mathcal{C}v_0 = (0)$ for any $v_0 \in V_{min} \cap V(+)$. But since $V^{\prime} = \{ v \in V \ | \ \mathcal{C}v = 0\}$ is an $E$-submodule of $V$, we are done by the irreducibility of $V$. \\  
    (2) By (1), $\mathcal{C}$ acts trivially on $V$. From Lemma \ref{L9.1}, we obtain $V(+) \neq (0)$. The uniform boundedness of $V(+)$ follows from the Weyl group argument presented above. 
\end{proof}

\brmk \label{R9.6}
The explicit realization of the highest weight space $V(+)$, which in turn determines $V$ uniquely, depends heavily on the graded subalgebra $\mathcal{D}$ of $\mathcal{S}_{n+1}$ used in the construction of $E$. For instance, if $E = \tau_{2}^{M}$ (see \S \ref{SS4.4}), then it can be deduced from \cite{E1} that $V(+) \cong \mathbb{C}[t_0^{\pm s_0},t_1^{\pm s_1}]$ for some $s_0, s_1 \in \mathbb{Z}$. But on the other hand, if we take $E = \tau_2^{\mathcal{S}}$ (see \S \ref{SS4.3}), then it immediately follows from \cite[Theorem 3.2]{CLT} that $V(+) \cong W \otimes \mathbb{C}[t_0^{\pm 1},t_1^{\pm 1}]$ for some finite-dimensional $\mathfrak{sl}_2$-module $W$.    
\ermk

\subsection{Restricted generalized highest weight modules} \label{SS9.2}
In this subsection, unless otherwise explicitly stated, $V$ will always denote an irreducible quasi-finite module over $E$ (or over $\tau_{n+1}^{tor}$ for $n>0$) with non-trivial action of $E_c$ (respectively $\overline{L}(\mathfrak{g})$).

\smallskip

\noindent \textbf{Notations.} For any $n \in \mathbb{Z_{+}}$, let $\{e_0, \ldots, e_n \}$ denote the canonical $\mathbb{Z}$-basis of $\mathbb{Z}^{n+1}$. If $\underline{k}, \underline{l} \in \mathbb{Z}^{n+1}$ such that $k_i  \geqslant l_i$ for all $0 \leqslant i \leqslant n$, then we shall say that $\underline{k} \geqslant \underline{l}$. For $p, q \in \mathbb{Z}$, let us set $[p,q] = \{m \in \mathbb{Z} \ | \ p \leqslant m \leqslant q \}$ and define $(-\infty, p]$ as well as $[q, \infty)$ similarly.
  
\blmma \label{L9.7}
For any $\underline{k} \in \mathbb{Z}^{n+1}, \ x_{\alpha} \otimes t^{\underline{k}} \in \mathfrak{g}_{\alpha} \otimes \mathbb{C}t^{\underline{k}} (\alpha \neq 0)$ either acts injectively or locally nilpotently on $V$.
\elmma

\begin{proof}
 Assume that there exists a non-zero $v \in V$ such that $(x_{\alpha} \otimes t^{\underline{k}})v = 0$ for some $\underline{k} \in \mathbb{Z}^{n+1}$ and $\alpha\in \Delta$. Since $\mathcal{D}$ acts on $L(\mathfrak{g})$ by derivations of $L(\mathfrak{g})$ and $\mathcal{C}$ commutes with $L(\mathfrak{g}) $, it is evident that for any $u_1, \ldots, u_r \in \mathcal{D}$ and $z_1, \ldots,z_r \in \mathcal{C}$, there exists $M(r) \in \mathbb{N}$ such that 
\begin{align*}
(x_{\alpha} \otimes t^{\underline{k}})^{M(r)}\big((u_1u_2u_3 \ldots u_r)v\big) = 0 \ \text{and} \  (x_{\alpha} \otimes t^{\underline{k}})\big((z_1z_2z_3 \ldots z_r)v\big) = 0.  
\end{align*} 
Morever as the adjoint action of $L(\mathfrak{g}) \oplus \mathcal{C}$ is integrable, it follows that for any $x_1, \ldots, x_s \in L(\mathfrak{g})$, there exists $M(s) \in \mathbb{N}$ satisfying $(x_{\alpha} \otimes t^{\underline{k}})^{M(s)}\big((x_1x_2x_3 \ldots x_s)v\big) = 0$. The lemma is now a direct consequence of the irreducibility of $V$. 
\end{proof}

\bdfn
$V$ is said to be a restricted generalized highest weight (GHW) module over $E$ if there exists a non-zero $v_{\Lambda_{0}} \in V_{\Lambda_{0}}$ such that $(E_{c})_{\underline{k}}v_{\Lambda_{0}} = (0) \ \forall \ \underline{k} \geqslant (N, \ldots,N)$ for some $N \in \mathbb{N}$, i.e. $V$ is a GHW module when restricted to $E_c$. In this case, $v_{\Lambda_0}$ is called a restricted GHW vector. 
\edfn

\bppsn \label{P9.9}
$V$ is either a level zero integrable module or a restricted GHW module.  
\eppsn

\begin{proof}
Suppose $V$ is not a level zero integrable module. Recall that $P_{D}(V) \subseteq \lambda + \mathbb{Z}^{n+1}$.   \\
\textbf{Claim.} There exists some $\underline{0} \neq \underline{r} \in \mathbb{Z}^{n+1}$ such that $\{\mathrm{dim}V_{\lambda + k \underline{r}} \ | \ k \in \mathbb{Z} \}$ is not uniformly bounded. \\
If not, then for any $\underline{0} \neq \underline{r} \in \mathbb{Z}^{n+1}$, $W(\underline{r}) = \oplus_{k \in \mathbb{Z}} V_{\lambda + k \underline{r}}$ is a uniformly bounded module over the affine Kac--Moody algebra, formed by considering the $\mathbb{Z}\underline{r}$-graded loop algebra along with its corresponding $1$-dimensional universal central extension and then accordingly adding the degree derivation. Then $W(\underline{r})$ is a uniformly bounded module and $x_{\alpha} \otimes t^{\underline{r}} \in \mathfrak{g}_{\alpha} \otimes \mathbb{C}t^{\underline{r}}$ acts locally nilpotently on $W(\underline{r})$ by Claim 2 of Lemma \ref{L8.3}. Now applying Lemma \ref{L7.8} and Lemma \ref{L9.7}, we can infer that $V$ is a level zero integrable module, which contradicts our initial assumption. Hence the claim. 

Assume that $\{\mathrm{dim}V_{\lambda - k \underline{r}} \ | \ k \in \mathbb{N} \}$ is not uniformly bounded. Put $\underline{r} = \underline{r_0}$ and extend it to a $\mathbb{Z}$-basis $\{\underline{r_i}\}_{i=0}^{n}$ of $\mathbb{Z}^{n+1}$. Put $N_j =$ max$\{|r_0^{j}|, |r_1^{j}|, |r_0^{j} + r_1^{j} + r_k^{j}| \ | \ 2 \leqslant k \leqslant n \}$ for each $0 \leqslant j \leqslant n$, where $\underline{r_i} = (r_i^{0}, \ldots, r_i^{n}) \ \forall \ 0 \leqslant i \leqslant n$. Take $N =$ max$\{N_j \ | \ 0 \leqslant j \leqslant n \} \in \mathbb{N}$ and consider 
\begin{align*}
B_N(\lambda) = \lambda + \{\underline{s} \in \mathbb{Z}^{n+1} \ | \ |s_i| \leqslant N \ \forall \ 0 \leqslant i \leqslant n \}.
\end{align*} 
Now since $\{\mathrm{dim}V_{\lambda - k \underline{r_0}} \ | \ k \in \mathbb{N} \}$ is not uniformly bounded, there exists some $k \in \mathbb{N}$ such that
\begin{align} \label{(7.1)}
\mathrm{dim}V_{\lambda - k\underline{r_0}} > M \sum_{\eta \in B_N(\lambda)} \mathrm{dim}V_{\eta}, \ \text{where} \ \mathrm{dim}  E_{\underline{r}} \leqslant M \ \forall \ \underline{r} \in \mathbb{Z}^{n+1}.
\end{align}
Set $\underline{s_0} = (k+1)\underline{r_0} + \underline{r_1}, \ \underline{s_1} = k\underline{r_0} + \underline{r_1}, \ \underline{s_j} = \underline{s_0} + \underline{r_j} \ \forall \ 2 \leqslant j \leqslant n$.
Then the linear map sending $\underline{r_i}$ to $\underline{s_i}$ has determinant $1$ and so $\{\underline{s_i} \}_{i=0}^{n}$ is a $\mathbb{Z}$-basis of $\mathbb{Z}^{n+1}$. Note that $\underline{s_i} + (\lambda - k\underline{r_0}) \in B_N(\lambda) \ \forall \ 0 \leqslant i \leqslant n$. As a result, it follows from (\ref{(7.1)}) that there exists $0 \neq v_0 \in V_{\lambda - k\underline{r_0}}$ such that $E_{\underline{s_i}}v_0 = 0 \ \forall \ 0 \leqslant i \leqslant n$. Hence up to a twist of an automorphism, $E_{e_i}v_0 = 0 \ \forall \ 0 \leqslant i \leqslant n$. Finally, applying Lemma \ref{L3.6} and the commutator relations on $E_c$, we can deduce that $V$ is a restricted GHW module.      
\end{proof}

\blmma\label{L9.10}
Let $V$ be an irreducible restricted GHW module over $E$ (or over $\tau_{n+1}^{tor}$ for $n>0$) having a restricted GHW vector $v_0$ with weight $\Lambda_0$. 
\begin{enumerate}
	\item Every non-zero $v \in V$ is a restricted GHW vector of $V$. 
	\item $(E_c)_{-\underline{k}}v \neq (0)$ for any $0 \neq v \in V$ and $\underline{k} \in \mathbb{N}^{n+1}$. 
	\item For each $\mu \in P_D(V)$ and $\underline{k} \in \mathbb{N}^{n+1}$, $\{m \in \mathbb{Z} \ | \ \mu + m \underline{k} \in P_D(V) \} = (-\infty, l]$ for some $l \in \mathbb{Z}_{+}$.
	\item (a) $\Lambda_{0} + \underline{k} \notin P_D(V) \ \forall \ \underline{0} \neq \underline{k} \in \mathbb{Z}_{+}^{n+1}$. \\ 
	(b) $\Lambda_{0} - \underline{k} \in P_D(V) \ \forall \   \underline{k} \in \mathbb{Z}_{+}^{n+1}$. \\ 
	(c) For any $\underline{k}, \underline{l} \in \mathbb{Z}^{n+1}$ with $\underline{k} \leqslant \underline{l}$, the condition $\Lambda_{0} + \underline{k} \notin P_D(V)$ implies $\Lambda_{0} + \underline{l} \notin P_D(V)$. \\ 
	(d) For any $\underline{0} \neq \underline{k} \in \mathbb{Z}_{+}^{n+1}$ and $\underline{l} \in \mathbb{Z}^{n+1}$, we have \\
	$\{m \in \mathbb{Z} \ | \ \Lambda_{0} +  \underline{l} + m \underline{k} \in P_D(V) \} = (-\infty, q]$ for some $q \in \mathbb{Z}$. 
\end{enumerate}
\elmma
\begin{proof}
(1) By hypothesis, there exists $N \in \mathbb{N}$ such that $(E_c)_{\underline{k}}v_{\Lambda_0} = (0) \ \forall \ \underline{k} \geqslant (N, \ldots,N)$. Let $\underline{m} \in \mathbb{Z}^{n+1}$ be arbitrary and set $N^{\prime} = N + |m_0| + |m_1| + \ldots + |m_n|$. Now as $E_c$ is an ideal of $E$,   
\begin{align} \label{(7.2)}
(E_c)_{\underline{k}} \big(E_{\underline{m}}v_{\Lambda_{0}} \big) \subseteq E_{\underline{m}} \big( (E_c)_{\underline{k}}v_{\Lambda_{0}} \big) + (E_c)_{\underline{k} + \underline{m}}v_{\Lambda_{0}} = (0) \ \forall \ \underline{k} \geqslant (N^{\prime}, \ldots, N^{\prime}).  
\end{align}
But since $V = U(E)v_{\Lambda_0}$ and $U(E)$ is generated by $\{E_{\underline{m}} \ | \ \underline{m} \in \mathbb{Z}^{n+1} \}$ (by the PBW theorem), the required result immediately follows from (\ref{(7.2)}). \\
(2) If possible, let $(E_c)_{-\underline{k}}v = 0$ for some $\underline{k} \in \mathbb{N}^{n+1}$ and $v \in V$. Then $v$ is a restricted GHW vector (by (1)) and so there exists $N \in \mathbb{N}$ such that $(E_c)_{\underline{m}}v = (0) \ \forall \ \underline{m} \geqslant (N, \ldots, N)$. This implies that $\big ((E_c)_{e_i + N \underline{k}} \big )v = 0 \ \forall \ 0 \leqslant i \leqslant n$. Now $\mathbb{Z}^{n+1}$ is generated by $\{e_i + N \underline{k}, - \underline{k} \ | \ 0 \leqslant i \leqslant n \}$ and thus $E_c$ is generated as a Lie algebra by $\{(E_c)_{-\underline{k}}, (E_c)_{e_i + N\underline{k}} \ | \ 0 \leqslant i \leqslant n \}$. This implies that $(E_c)v = 0$ and thus $E_c$ acts trivially on $V$ from Remark \ref{R7.3}, which is a contradiction. Hence the result follows. \\  
(3) Set $J = \{m \in \mathbb{Z} \ | \ \mu + m \underline{k} \in P_D(V) \}$. 
From (2), either $J = (-\infty, l]$ for some $l \in \mathbb{Z_{+}}$ or $J = \mathbb{Z}$. We claim that $J \neq \mathbb{Z}$. 
 If not, then $J = \mathbb{Z}$. Let $\mathfrak{L}_{\underline{k}}$ be the Lie algebra formed by identifying $\mathbb{Z}$ with $\mathbb{Z}\underline{k}$ and henceforth considering the corresponding $\mathbb{Z}\underline{k}$-graded (untwisted) loop algebra along with its $1$-dimensional universal central extension.
Then $\mathfrak{L}_{\underline{k}}$ is a finitely generated (affine) Lie algebra with finite-dimensional graded components. For any $R \in \mathbb{N}$, it can be verified that
\begin{align*}
[\mathfrak{L}_{\underline{k}}^{-}, (\mathfrak{L}_{\underline{k}})_{\geqslant R}] = \mathfrak{L}_{\underline{k}}, \ \text{where} \ \mathfrak{L}_{\underline{k}}^{\pm} = \bigoplus_{\pm p \in \mathbb{N}}(\mathfrak{L}_{\underline{k}})_p, \ (\mathfrak{L}_{\underline{k}})_{\geqslant R} = \bigoplus_{p \in \mathbb{Z}, \ p \geqslant R}(\mathfrak{L}_{\underline{k}})_p.   
\end{align*} 
Set $M_{\underline{k}} = \bigoplus_{p \in \mathbb{Z}} V_{\mu +
	p\underline{k}}$, which is a $\mathbb{Z}$-graded
$\mathfrak{L}_{\underline{k}}$-module. Again by (1), each $0 \neq v \in M_{\underline{k}}$ is a restricted GHW vector. So there exists $S \in \mathbb{N}$ such that $\big((\mathfrak{L}_{\underline{k}})_{\geqslant S} \big )v = (0)$, which implies that $\mathfrak{L}_{\underline{k}}^{+}v$ is finite-dimensional for each $v \in M_{\underline{k}}$. Now proceeding verbatim as in \cite[Lemma 1.6]{OM}, it can be deduced that for each $m \in \mathbb{Z}$, there exist $m^{\prime} > m$ and $0 \neq v(m^{\prime}) \in V_{\mu + m^{\prime}\underline{k}}$ such
that $\mathfrak{L}_{\underline{k}}^{+}v(m^{\prime}) = (0)$. But then the weight $\mu$ occurs in infinitely many irreducible highest weight $\mathfrak{L}_{\underline{k}}$-subquotients of $M_{\underline{k}}$ and hence $\dim V_{\mu} = \infty$, which establishes the claim. \\                
(4) (a), (b), (c) follow verbatim as in \cite[Lemma 3.6]{MZ}, while (d) can be deduced by considering the $\mathbb{Z}$-basis $\{e_i^{\prime}\}_{i=0}^{n}$ of $\mathbb{Z}^{n+1}$ as in \cite[Lemma 3.6]{MZ} and applying (3) with a change of co-ordinates.
\end{proof} 

\section{Quasi-finite modules over nullity 2 fgc EALAs} \label{S10} 
In this section, we classify the irreducible quasi-finite modules with non-trivial core action over an fgc EALA $E := E_2$ of nullity $2$ (or over the toroidal Lie algebra $\tau_2^{tor}$). Throughout this section, we shall always assume (using Lemma \ref{L7.4}) that $K_1$ acts trivially on our module.

\subsection{Highest weight type modules} \label{SS10.1}
In this subsection, we introduce the notion of highest weight type modules over $E_2$ and then classify their highest weight spaces.

Consider a $\mathbb{Z}$-grading on $E$ by means of eigenvalues of
$d_0$. Using Corollary \ref{C3.3} and Lemma \ref{L3.6}, this naturally gives rise to the following triangular decomposition of $E$.
\begin{align*} 
E = &\ E^{-} \ \oplus \ E^{0} \ \oplus \ E^{+}, \
\text{where}\\
E^{0} = &\ \sum_{r \in \Xi}\mathbb{C}t_1^rK_0 \oplus \sum_{r \in \Xi}\mathbb{C}t_1^{-r}d_0 \oplus \ \mathfrak{g} \otimes \mathbb{C}[t_1, t_1^{-1}] \ \oplus \ \mathbb{C}K_1 \ \oplus \ \mathbb{C}d_1, \ \text{for a subset} \ \Xi \ \text{of} \ \mathbb{Z} \ \text{containing} \ 0.
\end{align*}

Let $X$ be an irreducible $E^{0}$-module. Postulate $E^{+}X = (0)$, which in turn induces an $(E^{0} \oplus E^{+})$-module structure on $X$. Define the \textit{generalized Verma module} $M(X)$ over $E$  by setting
\begin{align} \label{Triangle}
M(X) = U(E) \otimes_{U(E^{0} \oplus E^{+})} X.
\end{align}     
By standard arguments, it can be shown that $M(X)$ has a unique irreducible quotient, say $L(X)$. Again there exists an action of $GL(2, \mathbb{Z})$ on $E$ (in the sense described in Section \ref{S6}). So for each $A \in GL(2, \mathbb{Z})$, we obtain a twisted irreducible module over $E$, which we shall denote by $L(X)^A$, where the action of $E$ is twisted by an automorphism $A$ of $\mathbb{Z}^{2}$. In this case, the twisted module will finally turn out to be a module over a Lie algebra isomorphic to $E$ (see \cite{BB,BZ} for examples of quasi-finite modules over $\mathbb{Z}^2$-graded Lie algebras). 

\brmk \label{R10.1}
In case of $\tau_2^{tor}$ (see \S \ref{SS4.2}), we also get a similar triangular decomposition given by $\tau_2^{tor} = (\tau_2^{tor})^{-} \oplus (\tau_2^{tor})^{0} \oplus (\tau_2^{tor})^{+}$, with $(\tau_2^{tor})^{0} = \sum_{r \in \mathbb{Z}} \mathbb{C}t_1^{r}K_0 \oplus \mathbb{C}d_0 \oplus \mathfrak{g} \otimes \mathbb{C}[t_1,t_1^{-1}] \oplus \mathbb{C}K_1 \oplus \mathbb{C}d_1$.  
\ermk
\bdfn \label{D10.2}
A linear map $z : X \longrightarrow X$ is said to be an $E^0$-central operator of degree $m$ if 
\begin{enumerate}
	\item $z$ commutes with the action of $(E^0)^{\prime}$, where \\  $(E^{0})^{\prime} := \sum_{r \in \Xi} \mathbb{C}t_1^rK_0 \oplus
	\sum_{r \in \Xi} \mathbb{C}t_1^{-r}d_0 \oplus \mathfrak{g} \otimes \mathbb{C}[t_1,t_1^{-1}] \oplus \mathbb{C}K_1$ and 
	\item $d_1 z -z d_1 = mz$.
\end{enumerate}
\edfn

\blmma \label{L10.3}
$t_1^{r}K_0$ and $t_1^{r}d_0$ are $E^{0}$-central operators of degree $r$ on $X$ for any $r \in \Xi$.
\elmma

\begin{proof}
It is enough to prove the assertion only for $t_1^rd_0$. First note that
\begin{align*}
\big (t_1^rd_0 \cdot (t_1^sK _0) \big )(t_1^pd_i) = -t_1^sK_0([t_1^rd_0, t_1^pd_1] \big) = ir \delta_{r+s+p,0} = r t_1^{r+s}K_1(t_1^pd_i) \ \forall \ i = 0, 1,
\end{align*} 
which implies that $t_1^rd_0 \cdot (t_1^sK_0) = rt_1^{r+s}K_1 = r \delta_{r+s,0}K_1$, due to Lemma \ref{L3.6}. Now the trivial action of $K_1$ gives $[t_1^rd_0,t_1^sK_0] = 0$. Moreover we can check that $\phi(t_1^rd_0,t_1^sd_0)(t_1^pd_i) = 0 \ \forall \ i=0,1$, which gives $[t_1^rd_0, t_1^sd_0] = 0$. Finally using (\ref{(3.1)}), we have $\phi(d_1, t_1^rd_0)(t_1^sd_i) = \phi(t_1^{r}d_0,t_1^sd_i)(d_1) = 0$, which yields $[d_1, t_1^{r}d_0] = rt_1^rd_0 + \phi(d_1,t_1^{r}d_0) = rt_1^rd_0$. Hence we are done as $[t_1^rd_0, x \otimes t_1^s] = 0 \ \forall \ x \in \mathfrak{g}, \ s \in \mathbb{Z}$.     
\end{proof}
\blmma \label{L10.4} \cite[Lemma 4.3, Lemma 4.4]{E2}
\
\begin{enumerate}
	\item If $z$ is an $E^0$-central operator of degree $m$ with $zx \neq 0$ for some non-zero $x \in X$, then $zx$ is injective.
	\item For a non-zero $E^0$-central operator $z$ of degree $m$ , there exists an $E^0$-central operator $T$ of degree $-m$ satisfying $Tz=zT=Id$.
\end{enumerate} 
\elmma


Recall that $\mathfrak{g} = \mathfrak{n}^{-} \oplus \mathfrak{h} \oplus \mathfrak{n}^{+}$ and consider another triangular decomposition of $E$ given by
\begin{align*} 
E = &\ E_{-} \ \oplus \ E_{0} \ \oplus \ E_{+}, \
\text{where}\\
E_{0} = &\ \sum_{r \in \Xi}\mathbb{C}t_1^rK_0 \oplus \sum_{r \in \Xi}\mathbb{C}t_1^{-r}d_0 \oplus \ \mathfrak{h} \otimes \mathbb{C}[t_1, t_1^{-1}] \ \oplus \ \mathbb{C}K_1 \ \oplus \ \mathbb{C}d_1 = E_0^{\prime} \oplus \mathbb{C}d_1.
\end{align*}
As in (\ref{Triangle}), we can now construct the \textit{generalized Verma module} with respect to this triangular decomposition and then obtain its unique irreducible $E$-quotient, which we shall denote by $V(T)$, where $T$ is an \textit{irreducible} $E_0$-module. Twisting this module by $B \in GL(2,\mathbb{Z})$ again gives rise to an irreducible module, say $V(T)^B$, over a Lie algebra isomorphic to $E$ (see Section \ref{S6}). We shall refer to these twisted modules of the form $V(T)^B$ as simply \textit{irreducible modules of the highest weight type}. Note that an analogous triangular decomposition also exists for $\tau_2^{tor}$, which similarly yields irreducible highest weight type modules over $\tau_2^{tor}$, induced by irreducible $(\tau_2^{tor})_0$-modules.

We now describe the irreducible quasi-finite modules over $E_0$ with trivial $K_1$-action. For any linear functional $\psi$ of $E_0^{\prime}$ with $\psi(K_1) = 0$, define an $E_0$-module structure on $\mathbb{C}[t_1,t_1^{-1}]$ by setting
\begin{align*}
(h \otimes t_1^r).t_1^s = \psi(h \otimes t_1^r)t_1^{r+s}, \,\, (t_1^pK_0).t_1^s = \psi(t_1^pK_0)t_1^{r+s}, \,\ (t_1^{-p}d_0).t_1^s = \psi(t_1^{-p}d_0)t_1^{-p+s}, \,\ d_1.t_1^s = st_1^s \
\end{align*}  
for all $h \in \mathfrak{h}, \ r, s \in \mathbb{Z}$ and $p \in \Xi$. Let $A_{\psi}$ be the $E_0$-submodule of $\mathbb{C}[t_1,t_1^{-1}]$ generated by $1$. For any $a \in \mathbb{C}$, let us now modify the action of $d_1$ on $A_{\psi}$ by simply taking the grade shift by $a$ and denote the resulting $E_0$-module by $A_{\psi,a}$. The following results can be readily deduced from \cite{C,E1}. 

\blmma \label{Lemma} 
\
\begin{enumerate}
\item For $\psi \in (E_0^{\prime})^*$ satisfying $\psi(K_1)=0$ and $a \in \mathbb{C}$, $A_{\psi,a}$ is irreducible over $E_0$ if and only if there exists $k \in \mathbb{N}$ such that $A_{\psi} = \mathbb{C}[t_1^{\pm k}]$ or $A_{\psi} = \mathbb{C}$ as vector spaces.  
\item Let $T$ be an irreducible quasi-finite module over $E_0$ such that $K_1$ acts trivially on $T$. Then $T \cong A_{\psi,a}$ for some $\psi \in (E_0^{\prime})^*$ with $\psi(K_1) = 0$ and  $a \in \mathbb{C}$.
\end{enumerate}
\elmma 

\brmk \label{Notation}
The above lemma conveys that $V(\psi,a) := V(A_{\psi,a})$ is an irreducible highest weight type module over $E$ (or $\tau_2^{tor}$) for any $\psi \in (E_0^{\prime})^*$ (respectively $((\tau_2^{tor})_0^{\prime})^*$) with $\psi(K_1) = 0$ and $a \in \mathbb{C}$.
\ermk
\subsection{Classification of irreducible restricted GHW modules} \label{SS10.2}
In this subsection, we show that every irreducible restricted GHW module $V$ over $E := E_2$ (or over $\tau_2^{tor}$), with a non-trivial action of $E_c$ (respectively $\overline{L}(\mathfrak{g})$), is a highest weight type module.

\blmma \label{L10.5}
Let $L(X)$ be an irreducible quasi-finite module over $E$, where $X$ is an irreducible $E^0$-module and $E_c$ acts non-trivially on $L(X)$.
If there exists $0 \neq m \in \Xi$ such that either $t_1^mK_0$ or $t_1^{-m}d_0$ acts non-trivially on $X$, then $X$ is a uniformly bounded $E^0$-module (with respect to $d_1$). 
\elmma
\begin{proof}
Put $\widehat{L}(\mathfrak{g}) = \mathfrak{g} \otimes \mathbb{C}[t_1,t_1^{-1}] \oplus \mathbb{C}K_1 \oplus \mathbb{C}d_1$. As $K_1$ acts trivially on $X$, both $t_1^rK_0$ and $t_1^{-r}d_0$ are $E^0$-central operators for all $r \in \Xi$. Subsequently every non-zero $t_1^rK_0$ and $t_1^{-r}d_0$ are invertible on $X$, with their inverses being $E^0$-central operators by Lemma \ref{L10.4}. \\
\textbf{Claim.} $X$ admits a (non-zero) finite-dimensional irreducible $(E^0)^{\prime}$-quotient. \\
Let $t_1^mK_0$ act non-trivially on $X$ for some $0 \neq m \in \Xi$. The case involving $d_0$ can be handled similarly. Now proceeding verbatim as in Claims 1, 2 and 3 of \cite[Theorem 4.5]{E2}, we can deduce that $X^{\prime} = \text{span} \{(t_1^mK_0)x - x \ | \ x \in X \}$ is a proper subspace of $X$ and $X/X^{\prime}$ is a finite-dimensional $(E^0)^{\prime}$-module. Also for any increasing sequence of $(E^0)^{\prime}$-submodules of $X$ containing $X^{\prime}$ given by
\begin{align*}
X^{\prime} \subsetneq X_1 \subsetneq \ldots \subsetneq \ldots,
\end{align*}
we get a decreasing chain of $(E^0)^{\prime}$-modules with $0 \neq \text{dim}(X/X_1) < \text{dim}(X/X^{\prime})$. But as $X/X^{\prime}$ is finite-dimensional, the above chain of  submodules terminates and so we get a maximal submodule of $X$ containing $X^{\prime}$, say $X^{\prime \prime}$. Then $\overline{X}:= X/X^{\prime \prime}$ is irreducible over $(E^0)^{\prime}$. Hence the claim.

For a fixed $\mu \in P_{D}(X)$, set $\beta = \mu(d_1) \in \mathbb{C}$ and put a $E^0$-module structure on $\overline{X} \otimes \mathbb{C}[t_1, t_1^{-1}]$ as follows:
\begin{align} \label{(8.1)}
(y \otimes t_1^r)(\overline{x} \otimes t_1^s) = (y \otimes t_1^r)\overline{x} \otimes t_1^{r+s},  &\ &\ d_1(\overline{x} \otimes t_1^s) = (s + \beta)(\overline{x} \otimes t_1^s), \\
t_1^pK_0(\overline{x} \otimes t_1^s) = (t_1^pK_0)\overline{x} \otimes t_1^{p+s}, &\ &\ t_1^{-p}d_0(\overline{x} \otimes t_1^s) = (t_1^{-p}d_0)\overline{x} \otimes t_1^{s-p}. 
\end{align}
For any $\overline{x} \in \overline{X}$ and $s \in \mathbb{Z}$, put $\overline{x}(r) = \overline{x} \otimes t^r$.
The action of $d_1$ induces a $\mathbb{Z}$-grading on $X$, say  $X = \oplus_{r \in \mathbb{Z}} X_r$. Thus we have a natural map
\begin{equation} \label{(4.1)} 
\phi : X \longrightarrow \overline{X} \otimes \mathbb{C}[t_1,t_1^{-1}]
\end{equation}
$\,\,\,\,\,\,\,\,\,\,\,\,\,\,\,\,\,\,\,\,\,\,\,\,\,\,\,\,\,\,\,\,\,\,\,\,\,\,\,\,\,\,\,\,\,\,\,\,\,\,\,\,\,\,\,\,\,\,\,\,\,\,\,\,\,\,\,\,\,\,\,\,\,\,\,\,\,\,\,\,\,\,\,\,\,\,\,\,\,\,\,\,\,\,\,\,\,\,\,\,\,\,\,\ x \longmapsto \ \overline{x} \otimes t_1^{r},\ x \in X_{r}$. \\
Then $\phi$ is a non-zero $E^0$-module homomorphism and so $\phi$ must be injective, due to the irreduciblity of $X$. This implies that $X \cong \phi(X)$ is an irreducible $E^0$-submodule of $\overline{X} \otimes \mathbb{C}[t_1,t_1^{-1}]$, which is clearly uniformly bounded with respect to $d_1$, with $d_0$ acting by a fixed scalar on $X$ (as $[d_0,E^0]=0$).
 \end{proof}

\bppsn \label{P10.11}
Let $\lambda_0 \in P_D(V)$ such that $(\lambda_0 + \mathbb{N}e_0 +
\mathbb{Z}e_1) \bigcap P_D(V) = \emptyset$. Then $V \cong L(X)^A$ for an irreducible uniformly bounded module $X$ over $E^0$ and  some $A \in GL(2, \mathbb{Z})$. 
\eppsn
\begin{proof}
	By hypothesis, $V_{+} := \{v \in V \ | \ E^{+}v = (0) \} \neq (0)$. Due to Lemma \ref{L7.4}, we can assume that $K_1$ acts trivially on $V$. Now since $V$ is irreducible, an application of the PBW theorem implies that $V_{+}$ is an irreducible quasi-finite module with respect to $\mathbb{C}d_1$ over $E^0$. Also note that
	\begin{align*}
		V_{+} = \bigoplus_{m \in \mathbb{Z}} (V_{+})_{m}, \ \text{where} \ (V_{+})_{m} = \{v \in V_{+}\ | \ d_1v = (\lambda_0(d_1)+m)v \} \ \forall \ m \in \mathbb{Z}.
	\end{align*}
	\textbf{Claim 1.} $V_{+}$ is either uniformly bounded or a highest weight module or a lowest weight module. \\
	By Lemma \ref{L10.5}, $V_{+}$ is uniformly bounded over $E^0$ if there exists some $0 \neq m \in \Xi$ such that either $t_1^mK_0$ or $t_1^{-m}d_0$ acts non-trivially on $V_{+}$. Again if both $t_1^mK_0$ and $t_1^{-m}d_0$ act trivially on $V_{+}$ for each  $0 \neq m \in \Xi$, then $V_{+}$ is irreducible over $\widetilde{L}(\mathfrak{g}):= \mathfrak{g} \otimes \mathbb{C}[t_1,t_1^{-1}] \oplus \mathbb{C}d_1$, with $K_0$ and $d_0$ acting by fixed scalars on $V_{+}$ and so the claim follows from our Claim 1 of Lemma \ref{L8.3}. \\
    \textbf{Claim 2.} $V_{+}$ is a uniformly bounded module over $E^0$.\\
	If not, then  $V_{+}$ is either a highest weight module or a lowest weight module by Claim 1. Without loss of generality, let $V_{+}$ be a highest weight module. So there exists some non-zero $v_0 \in V_{+}$ such that $(E^{0})_{j}v_{0} = (0) \ \forall \ j \in \mathbb{N}$. Moreover by Lemma \ref{L10.5}, both $t_1^mK_0$ and $t_1^{-m}d_0$ must act trivially on $V_{+} \ \forall \ m \in \Xi \setminus \{0\}$. Consequently $V_{+}$ is also an irreducible highest weight module over $\widetilde{L}(\mathfrak{g})$, as $K_0$ and $d_0$ act by fixed scalars on $V_{+}$.  Now if $\mathfrak{g} \otimes \mathbb{C}[t_1, t_1^{-1}]$ acts trivially on $V_{+}$, then $V_{+}$ is a $1$-dimensional irreducible module over the abelian Lie algebra $\mathbb{C}K_0 \oplus \sum_{i=0}^{1} \mathbb{C}d_i$ and thus we are done by considering $X=V_{+}$. So let us suppose $\mathfrak{g} \otimes \mathbb{C}[t_1, t_1^{-1}]$ acts non-trivially on $V_{+}$. Then there exist a vector $v_0$, say of weight $\mu_0$, in (the highest weight space of) $V_{+}$ and $h_0 \in \mathfrak{h}$ such that $h_0$ acts by a non-zero scalar on $v_0$. As a result, $V = U(E)v_0$ is a highest weight module with
	\begin{align}\label{P(V)}
	 P_D(V) \subseteq(\mu_0 - \mathbb{N}e_0 + \mathbb{Z}e_1) \bigcup (\mu_0 - \mathbb{Z_{+}}e_1).
	 \end{align}
	 Therefore $V$ is a restricted GHW module with GHW vector $v_0$ satisfying
	 \begin{align} \label{(5.1)}
	 (\mu_{0} + \mathbb{N}e_0 + \mathbb{Z}e_1) \bigcap
	 P_D(V) = \emptyset, \,\,\,\ (\mu_{0} + \mathbb{Z}e_1) \bigcap
	 P_D(V) \neq \emptyset.
	 \end{align}
	  We now show that the $\widetilde{L}(\mathfrak{g})$-module
	 $V^{\prime} := \bigoplus_{m \in \mathbb{Z}} V_{\mu_0 +
	 me_1}$ is uniformly bounded. \\
 Pick any non-zero $v \in V_{\mu_0 + m^{\prime}e_1}$, where $m^{\prime} \in \mathbb{Z}$ is arbitrary. If possible, let $(E_c)_{-m^{\prime}e_1 - e_0}v =(0)$. \\ Now
 $\mathbb{Z}^{2}$ is generated by the elements $\{me_1
 + e_0,  -m^{\prime}e_1 - e_0 \ | \ m \in
 \mathbb{Z} \}$ for any $m^{\prime} \in
 \mathbb{Z}$. Subsequently $E_c$ is generated as a Lie
 algebra by $\{(E_c)_{-m^{\prime}e_1 - e_0}, \ (E_c)_{me_1 + e_0} \ | \ m \in
 \mathbb{Z} \}$. In view of (\ref{(5.1)}), $E_c$ acts trivially on $v$. But then by Remark \ref{R7.3}, $E_c$ acts trivially on
 $V$, which contradicts our initial assumption. This shows that  
 $(E_c)_{-m^{\prime}e_1 - e_0}v \neq (0)$. \\
Thus for every $m^{\prime} \in \mathbb{Z}$,
 we have an injective linear map
 \begin{align*}
 T_{m^{\prime}}: V_{\mu_0 +
 	m^{\prime}e_1} \longrightarrow (V_{\mu_0 -
 	e_0})^{\oplus N}, \ \text{where} \ N = \dim (E_c)_{-m^{\prime}e_1 - e_0},
\end{align*}     
which implies that $V^{\prime}$ is uniformly bounded. Clearly $V^{\prime}$ contains the non-trivial $\widetilde{L}(\mathfrak{g})$-submodule $V^{\prime \prime} = U(\widetilde{L}(\mathfrak{g}))v_{0}$. Using Zorn's lemma, we can obtain a non-zero $\widetilde{L}(\mathfrak{g})$-irreducible quotient of $V^{\prime \prime}$, say $W^{\prime \prime}$. Again since $h_0 \in \mathfrak{h}$ acts by a non-zero scalar at $v_{0}$, $W^{\prime \prime}$ is also a uniformly bounded $\widetilde{L}(\mathfrak{g})$-module, with $
\mathfrak{g} \otimes \mathbb{C}[t_1,t_1^{-1}]$ acting non-trivially on $W^{\prime \prime}$. This implies that $\mu_0$ is not a highest weight of $W^{\prime \prime}$. Thus there exists $k \in \mathbb{N}$ such that $\mu_0 + ke_1 \in 
P_D(W^{\prime \prime}) \subseteq P_D(V)$, which contradicts (\ref{P(V)}). This contradiction establishes our claim and also proves the proposition by taking $X = V_{+}$.    
\end{proof}

\blmma \label{L10.12}
In each of the following cases, $V \cong L(X)^A$ for an irreducible uniformly bounded module $X$ over $E^0$ and some $A \in GL(2, \mathbb{Z})$.
\begin{enumerate}
	\item There exists some $(p_0,p_1) \in \mathbb{N}^2$ and $\underline{m}, \underline{k} \in \mathbb{Z}^2$, where $k_0, k_1$ are relatively prime, such that
	\begin{align*}
	\bigg \{\Lambda_0 + \sum_{i=0}^{1}m_ie_i + \sum_{i=0}^{1}r_ip_ie_i \ | \
	\underline{r} \in \mathbb{Z}^{2}, \ \sum_{i=0}^{1}k_ip_ir_i = 0 \bigg
	\} \bigcap P_D(V) = \emptyset.
	\end{align*}
	\item There exist $(i,j) \in \mathbb{Z}^2$ and $(k,l) \in \mathbb{Z}^2 \setminus \{(0,0)\}$ such that
	\begin{align*}
	\{\Lambda_0 + (i,j) + m(k,l) \ | \ m \in \mathbb{Z} \}
	\bigcap P_D(V) = \emptyset.
	\end{align*}
	\item There exist $(i,j) \in \mathbb{Z}^2, \ (k,l) \in \mathbb{Z}^2 \setminus \{(0,0)\}$ and $p, q \in \mathbb{Z}$ such that 
	\begin{align*}
		\{m \in \mathbb{Z} \ | \ \Lambda_{0} + (i,j) + m(k,l) \in P_D(V) \} \supseteq (-\infty, p] \cup [q, \infty).
	\end{align*}
	
	\item There exist $(i,j), (k,l) \in \mathbb{Z}^2$ and $m_1, m_2, m_3 \in \mathbb{Z}$ with $m_1 < m_2 < m_3$ such that 
	\begin{align*}
		& \Lambda_0 + (i,j) + m_1 (k,l) \notin P_D(V), \\
		& \Lambda_0 + (i,j) + m_2(k,l) \in P_D(V), \\
		& \Lambda_0 + (i,j) + m_3(k,l) \notin P_D(V).
	\end{align*}
\end{enumerate}
\elmma

\begin{proof}
(1) Using Lemma \ref{L9.10}, we can show that, up to a twist of an automorphism, there exists a
unique $m_0 \in \mathbb{Z}$ satisfying the following properties
(see (1) and (2) of \cite[Lemma 3.3]{LZ1}). 
\begin{align*}
\big \{\Lambda_{0} + \sum_{i=0}^{1}l_ie_i \in P_D(V) \ | \ \underline{l} \in \mathbb{Z}^{2}, \ \sum_{i=0}^{1}k_il_i \geqslant m_0 \big \} = \emptyset \ \text{and} \\
S := \big \{\Lambda_{0} + \sum_{i=0}^{1}l_ie_i \in P_D(V) \ | \ \underline{l} \in \mathbb{Z}^{2}, \ \sum_{i=0}^{1}k_il_i = m_0-1 \big \} \neq \emptyset.
\end{align*}
Now since $k_0, k_1$ are relatively prime, there exist $s_0, s_1 \in \mathbb{Z}$ such that $\sum_{i=0}^{1}k_is_i = 1$. Taking $e_0^{\prime} = (s_0,s_1)$ and $e_1^{\prime} = (k_1,-k_0)$ yields a $\mathbb{Z}$-basis $\{e_0^{\prime}, e_1^{\prime} \}$ of $\mathbb{Z}^2$. Then for any $\Lambda_0 \in S$, we have $P_D(V) \cap (\Lambda_0 + \mathbb{N}e_0^{\prime} + \mathbb{Z}e_1^{\prime}) = \emptyset$ and hence the required result follows from Proposition \ref{P10.11}. \\
(2), (3) and (4) can be deduced using Lemma \ref{L9.10}, Proposition \ref{P10.11} and (1) as in \cite{LZ1,SY}.   		
\end{proof}

By Lemma \ref{L10.12}, it follows that if $V \ncong L(X)^A$ for any irreducible uniformly bounded $E^0$-module $X$ and $A \in GL(2, \mathbb{Z})$, then for any $(i,j) \in \mathbb{Z}^2$ and $(k,l) \in \mathbb{Z}^2 \setminus \{(0,0)\}$, there exists $p \in \mathbb{Z}$ such that
\begin{align}\label{4.3}
\{m \in \mathbb{Z} \ | \ \Lambda_{0} + (i,j) + m(k,l) \in P_D(V) \} = (-\infty, p] \ \text{or} \ [p, \infty).
\end{align}
Consequently from Lemma \ref{L9.10}, it follows that for each $i \in \mathbb{N}$, there exist $a_i, b_i \in \mathbb{Z}_{+}$ satisfying

\begin{align*}
b_i = &\ \max \{b \in \mathbb{Z} \ | \ \Lambda_0 + (-i,b) \in P_D(V) \}, \,\,\,\	a_i = \ \text{max} \{a \in \mathbb{Z} \ | \ \Lambda_0 + (a,-i) \in P_D(V) \}.
\end{align*} 
Then we can deduce the following results from \cite[Claim 1]{LZ1} and \cite[Claim 2]{LZ1} (also see \cite{SY}).\\
(R1) The following limits exist finitely and are finite.
\begin{align*}
\alpha = \lim_{k \rightarrow \infty} \dfrac{b_k}{k}, \quad
\beta = \lim_{k \rightarrow \infty} \dfrac{a_k}{k}.	
\end{align*}
(R2) $\alpha = \beta^{-1}$ is a positive irrational number.	\\
(R3) Define a total order $>_{\alpha}$ on $\mathbb{Z}^2$ by setting
\begin{align*}
(i,j) >_{\alpha} (k,l) \iff i \alpha + j > k \alpha + l.	
\end{align*}
This order on $\mathbb{Z}^2$ is dense, which means that for every $(k,l) >_{\alpha} (0,0)$, there exist infinitely many $(i,j) \in \mathbb{Z}^2$ such that $(0,0) <_{\alpha} (i,j) <_{\alpha} (k,l)$.\\
Let us put $\mathbb{Z}^2(+) = \{(i,j) \in \mathbb{Z}^2 \ | \ (i,j) >_{\alpha} (0,0) \}$ and 
$\mathbb{Z}^2(-) = \{(i,j) \in \mathbb{Z}^2 \ | \ (i,j) <_{\alpha} (0,0) \}$.\\ 
(R4) $ \Lambda_{0} + (i,j) \in P_D(V) \implies \Lambda_{0} + (k,l) \in P_D(V) \ \forall \ (k,l) <_{\alpha} (i,j)$.\\
The ordering $>_{\alpha}$ naturally induces a triangular decomposition of $E$, say $E = E_{>_{\alpha}}^{-} \oplus E_{>_{\alpha}}^{0} \oplus E_{>_{\alpha}}^{+}$, where $E_{>_{\alpha}}^{0} = \mathfrak{g} \oplus \sum_{i=0}^{1} \mathbb{C}K_i \oplus \sum_{i=0}^{1} \mathbb{C}d_i$. 

\blmma \label{L10.13}
$(\mathfrak{g} \otimes t^{-\underline{s}})v_{\mu} \neq (0)$ for any $\underline{s} = (a,b) \in \mathbb{Z}^2(+)$ and $0 \neq v_{\mu} \in V_{\mu}$, where $\mu \in P_D(V)$.
\elmma
\begin{proof}
Let $V \ncong L(X)^A$  for any irreducible uniformly bounded $E^0$-module $X$ and $A \in GL(2, \mathbb{Z})$. From (R3), (R4) and (\ref{4.3}), it then follows that for any $(i,j) \in \mathbb{Z}^2$, there exists $p \in \mathbb{Z}$ such that
	\begin{align}\label{(4.10)}
		\{m \in \mathbb{Z} \ | \ \Lambda_{0} + (i,j) + m(k,l) \in P_D(V) \} = (-\infty, p] \ \forall \ (k,l) \in \mathbb{Z}^2(+).
	\end{align}
	This implies that for any fixed $\underline{s} = (a,b) \in \mathbb{Z}^2(+)$, we have $(\mathfrak{g} \otimes t^{r\underline{s}})v_{\mu} = (0)$ for $r \in
	\mathbb{N}$ large enough. \\
	Put $L(\mathfrak{g}) = \mathfrak{g} \otimes \mathbb{C}[t_0^{\pm 1}, t_1^{\pm 1}]$ and let, if possible, $L(\mathfrak{g})_{-\underline{s}}v_{\mu} = (0)$ for some non-zero vector $v_{\mu} \in V_{\mu}$. 
	Taking $c = \gcd(a,b) \in \mathbb{N}$, we have $\underline{s} = c(a^{\prime}, b^{\prime})$ where $\gcd(a^{\prime},b^{\prime}) = 1$. So
	there exist $p, q \in \mathbb{Z}$ such that $a^{\prime}q - b^{\prime}p = 1$. Set $e_0^{\prime} = (a^{\prime},b^{\prime})$ and $e_1^{\prime} = (p,q)$, whence $\{e_0^{\prime}, \ e_1^{\prime}\}$ forms a $\mathbb{Z}$-basis of $\mathbb{Z}^2$. Now we can directly appeal to (\ref{4.3}) to infer that, for any $0 \neq l \in \mathbb{Z}$, there exists $r_l \in \mathbb{Z}$ such that 
	\begin{align*}
	T_l := \{k \in \mathbb{Z} \ | \ \mu + le_1^{\prime} + ke_0^{\prime} \in P_D(V) \} = (-\infty, r_l] \ \text{or} \ [r_l, \infty).
	\end{align*}  
	It suffices to only consider $T_l = (-\infty, r_l]$, as a similar argument will also work for the other case. This gives 
	$L(\mathfrak{g})_{le_1^{\prime} + cs_le_0^{\prime} \pm e_0^{\prime}}v_{\mu} =$ for large enough $s_l \in \mathbb{N}$. Now as $\mathrm{dim} \ \mathcal{C}_{\underline{m}} \leqslant 1 \ \forall \ \underline{m} \in \mathbb{Z}^2 \setminus \{(0,0) \}$ (by Lemma \ref{L3.6}) and $L(\mathfrak{g})_{-ce_0^{\prime}}v_{\mu} = 0$, we can use the commutator relations on $E_c$ to show that  $L(\mathfrak{g})_{le_1^{\prime} \pm e_0^{\prime}}v_{\mu} = 0$ for any $0 \neq l \in \mathbb{Z}$. As a result, we have $L(\mathfrak{g})_{\pm (e_0^{\prime} + e_1^{\prime})}v_{\mu} = L(\mathfrak{g})_{\pm (e_0^{\prime} + 2e_1^{\prime})}v_{\mu} = 0$. On the other hand, since $\mathcal{C} = \big [\mathfrak{h} \otimes \mathbb{C}[t_0^{\pm 1}, t_1^{\pm 1}],\mathfrak{h} \otimes \mathbb{C}[t_0^{\pm 1}, t_1^{\pm 1}] \big ]$ and  $\{e_0^{\prime} + e_1^{\prime}, e_0^{\prime} + 2e_1^{\prime}\}$ is a $\mathbb{Z}$-basis of $\mathbb{Z}^2$, it can be deduced that $(E_c)v_{\mu} = 0$. The irreducibility of $V$ then implies that $E_c$ acts trivially on $V$ (by Remark \ref{R7.3}), which is a contradiction and hence the lemma is proved.   
\end{proof}

\blmma \label{L10.14}
$\big (\mu + \mathbb{Z}^{2}(+) \big) \bigcap P_D(V) \neq \emptyset \ \forall \ \mu \in P_D(V)$.
\elmma
\begin{proof}
	Let us assume the contrary. Then there exists $\mu = \lambda + \delta_{\underline{s}} \in P_D(V)$ for some $\underline{s} \in \mathbb{Z}^2$, with $\big (\mu + \mathbb{Z}^{2}(+) \big) \bigcap P_D(V) = \emptyset$ and thus $V$ is a highest weight module relative to the triangular decomposition in (R4) with its highest weight space $V_{+}^{> \alpha} = \{v \in V \ | \ E_{>_{\alpha}}^{+}v = 0 \} = (V_{+}^{> \alpha})_{\mu} \neq (0)$. Using the PBW theorem and the irreducibility of $V$, we can deduce that $V_{+}^{> \alpha}$ is irreducible over $E_{>_{\alpha}}^{0}$. Note that dim$V_{+}^{> \alpha} < \infty$, as $V$ is quasi-finite. Again since $V_{\mu}$ is a finite-dimensional $\mathfrak{g}$-module, it is a standard fact that $V_{\mu}$ has a minimal weight, say $\eta \in P_{\mathfrak{h}}(V_{\mu})$ (see Definition \ref{D2.2}). For example, every finite-dimensional $\mathfrak{sl}_2$-module has either $0$ or $1$ as one of its weights, which are precisely the minimal weights of $\mathfrak{sl}_2$. As a result, $\eta(\gamma^{\vee}) \in \{0,1 \}$ for any $\gamma \in \Delta_{\mathfrak{g}}^{+}$, where $\text{span} \{x_{\gamma}, y_{\gamma}, \gamma^{\vee} \} \cong \mathfrak{sl}_2$. Fix any $0 \neq v_0 \in V_{\eta + \delta_{\underline{s}}}$. \\ 
	Now for any $k \in \mathbb{N}$, pick $\underline{m} = (a,b) \in \mathbb{Z}^2$ with $\dfrac{-1}{4k} <a \alpha + b < 0$. Due to (R3), there exists
	\begin{align}\label{(6.6)}
    \text{infinitely many} \ a \in \mathbb{Z} \ \text{and infinitely many} \ b \in \mathbb{Z} \ \text{with} \ (a,b) \ \text{satisfying the above relation}.
	\end{align}
	 Put $(a_i,b_i) = (2i-1)\underline{m} \ \forall \ 1 \leqslant i \leqslant k$. We can check that $(0,-1) = -e_1 <_{\alpha} (a_i,b_i) <_{\alpha} (0,0)$. Again set $\underline{p_i} = (1-2i)\underline{m} >_{\alpha} (0,0)$ and  $\underline{q_i} = e_1 - \underline{p_i} >_{\alpha} (0,0)$ for each $1 \leqslant i \leqslant k$. Observe that    
	\begin{align}\label{4.4}
		\underline{p_j} + \underline{q_j} = e_1 \ , \ \underline{q_1} >_{\alpha} \underline{q_i} \ , \ \underline{q_1} >_{\alpha}   \underline{p_j} \ \forall \ 2 \leqslant i \leqslant k, \ 1 \leqslant j \leqslant k, \\ 	
	\label{4.6}
		e_1 >_{\alpha} (0,0) , \	[E_c,
		\mathcal{C}] = (0) \ \text{and} \ \big (\mu +
		\mathbb{Z}^{2}(+) \big) \bigcap P_D(V) = \emptyset.
	\end{align}  
	\textbf{Claim 1.} $K_0, \ K_1$ acts trivially on $V_{+}^{> \alpha}$. \\
	By Lemma \ref{L7.4}, it is enough to show that $K_0$ acts trivially on $V_{+}^{> \alpha}$. If not, then $K_0$ acts by a non-zero scalar $c_0$ on $V$.
	For any $k \in \mathbb{N}$,
	consider $\{(x_{\gamma} \otimes t^{-\underline{q_i}})(y_{\gamma} \otimes t^{-\underline{p_i}}
	)v_0 \}_{1 \leqslant i \leqslant k}$, all of which belong to
	$V_{\mu - e_1}$. Then these vectors cannot be linearly
	independent for each $k \in \mathbb{N}$, else we shall have $\dim V_{\mu -e_1} = \infty$. So there exists $k^{\prime} \in \mathbb{N}$ satisfying
	$\sum_{i=1}^{k^{\prime}}\beta_i(x_{\gamma} \otimes t^{-\underline{q_i}})(y_{\gamma} \otimes t^{-\underline{p_i}})v_0 = 0$ for some $\underline{0} \neq \underline{\beta} = (\beta_1, \ldots, \beta_{k^{\prime}}) \in \mathbb{C}^{k^{\prime}}$. Take $\underline{r} = e_1 + 2\underline{m} >_{\alpha} \underline{0}$ and note that 
	\begin{align}\label{4.5}
	\underline{q_1} >_{\alpha} \underline{r} \ , \ \underline{q_i} <_{\alpha} \underline{r} \ , \ \underline{p_j} <_{\alpha} \underline{r} \ , \ e_1 >_{\alpha} \underline{r} \ \forall \ 2 \leqslant i \leqslant k^{\prime}, \ 1 \leqslant j \leqslant k^{\prime}.
	\end{align}
	An application of $\gamma^{\vee} \otimes t^{\underline{r}}$ to the above equation together with (\ref{4.4}), (\ref{4.6}), (\ref{4.5}) yields that 
	\begin{align} \label{E1}
    \beta_1[2(x_{\gamma} \otimes t^{\underline{r} - \underline{q_1}})(y_{\gamma} \otimes t^{- \underline{p_1}}) + a_1 \big (\sum_{i=0}^{1} t^{\underline{r} - \underline{q_1}}K_i \big )(y_{\gamma} \otimes t^{- \underline{p_1}})]v_0 \\ +  \beta_1^{\prime}[2\gamma^{\vee} \otimes t^{\underline{r}-e_1} + a_2 (t^{\underline{r} -e_1}K_0) +
    a_3 (t^{\underline{r} -e_1}K_1)]v_0 = 0, 
	\end{align}
	where $\beta_1^{\prime} = \sum_{i=2}^{k^{\prime}} \beta_i$ and some $a_1, a_2, a_3 \in \mathbb{C}$. Consequently applying $\gamma^{\vee} \otimes t^{e_1 - \underline{r}}$ to the above equation and using (\ref{4.4}), (\ref{4.6}) and (\ref{4.5}), we obtain
	\begin{align} \label{E2}
     \big (\eta(\gamma^{\vee}) - ac_0(x_{\gamma}|y_{\gamma}) \big) \beta_1 - \big (ac_0(\gamma^{\vee}|\gamma^{\vee}) \big ) \beta_1^{\prime} = 0. 
	\end{align}
    Again apply $(y_{\gamma} \otimes t^{\underline{p_1}})(x_{\gamma}\otimes t^{\underline{r} - \underline{q_1}})$ to (\ref{E1}), (10.13) and use (\ref{4.4}), (\ref{4.6}), (\ref{4.5}) along with $(x_{\gamma}|x_{\gamma}) = 0 = (y_{\gamma}|y_{\gamma})$ to get
    \begin{align} \label{E3}
    \big(2 \eta(\gamma^{\vee}) - 2ac_0(x_{\gamma}|y_{\gamma}) - \eta(\gamma^{\vee})^2 + a^2c_0^2(x_{\gamma}|y_{\gamma})^2 \big) \beta_1 + \big(2 \eta(\gamma^{\vee}) - 2ac_0(x_{\gamma}|y_{\gamma}) \big) \beta_1^{\prime} = 0. 
    \end{align}
    \noindent \textit{Case 1.} $\eta(\gamma^{\vee}) = 1$. \\
    From (\ref{E2}) and (\ref{E3}), we get 
    \begin{align} \label{E4}
    \big (1 - ac_0 \big (x_{\gamma}|y_{\gamma}) \big ) \beta_1 - \big (ac_0(\gamma^{\vee}|\gamma^{\vee}) \big )\beta_1^{\prime} = 0 = \big (1 - ac_0(x_{\gamma}|y_{\gamma}) \big)^2 \beta_1 + \big (2 - 2ac_0(x_{\gamma}|y_{\gamma}) \big ) \beta_1^{\prime}.
    \end{align}
    \textit{Case 2.} $\eta(\gamma^{\vee}) = 0$. \\
    Using (\ref{(6.6)}), let us choose $a \neq 0$. Again as $c_0(x_{\gamma}|y_{\gamma}) \neq 0$, it follows from (\ref{E2}) and (\ref{E3}) that
    \begin{align} \label{E5}
    \big (2 - ac_0(x_{\gamma}|y_{\gamma}) \big ) \beta_1 + 2 \beta_1^{\prime} = 0 = (x_{\gamma} | y_{\gamma}) \beta_1 + (\gamma^{\vee} | \gamma^{\vee}) \beta_1^{\prime}. 
    \end{align}
    Finally, since $c_0(x_{\gamma}|y_{\gamma})(\gamma^{\vee}|\gamma^{\vee}) \neq 0$, we can invoke (\ref{(6.6)}) to pick some $a \in \mathbb{Z}$ such that the system of linear equations in both (\ref{E4}) and (\ref{E5}) admit a unique solution. This implies that $\beta_1 = \beta_1^{\prime} = 0$. Similarly we can deduce that $\beta_2 = \ldots = \beta_{k^{\prime}} = 0$, which is a contradiction. Hence the claim. \\
 \noindent \textbf{Claim 2.} $\mathfrak{g}$ acts trivially on $V_{+}^{> \alpha}$. \\
If not, then there exist $0 \neq v_0 \in V_{+}^{> \alpha}$ and $\lambda^{\prime} \in P_{\mathfrak{g}}^{+}$ such that $\gamma^{\vee}v_0 = \lambda^{\prime}(\gamma^{\vee})v_0 \neq 0$ for some $\gamma^{\vee} \in \mathfrak{h}, \ \gamma \in \Delta_{\mathfrak{g}}^{+}$. Then as in the last claim, there exist $k^{\prime} \in \mathbb{N}$ and $\underline{0} \neq (\alpha_1, \ldots, \alpha_{k^{\prime}}) \in \mathbb{C}^{k^{\prime}}$ with $\sum_{i=1}^{k^{\prime}}\alpha_i(x_{\gamma} \otimes t^{-\underline{q_i}})(y_{\gamma} \otimes t^{-\underline{p_i}}
)v_0 = 0$. Applying $(\gamma^{\vee} \otimes t^{e_1 - \underline{r}})(\gamma^{\vee} \otimes t^{\underline{r}})$ to this equation and using (\ref{4.4}), (\ref{4.6}) and (\ref{4.5}), it can be deduced using Claim 1 (see (\ref{E2})) that 
$\alpha_1 \lambda^{\prime}(\gamma^{\vee}) = 0$, which gives $\alpha_1 = 0$. Similarly we can show that $\alpha_2 = \ldots = \alpha_{k^{\prime}} = 0$. This contradiction proves the claim. 


Consequently it follows from Claim $1$ and Claim $2$ that $E_c$ acts trivially on $V$, as $V$ is an irreducible highest weight module relative to the triangular decomposition in (R4). This contradicts our hypothesis, which thereby proves the lemma.

\end{proof}

\bthm \label{T10.15}
\
\begin{enumerate}
\item $V \cong L(X)^A$ for a uniformly bounded irreducible $E^0$-module $X$ and $A \in GL(2, \mathbb{Z})$.
\item $V \cong V(\psi,a)^B$ for some $\psi \in (E_0^{\prime})^*$ (or $((\tau_2^{tor})_0^{\prime})^*$) with $\psi(K_1) = 0, \ a \in \mathbb{C}$ and $B \in GL(2, \mathbb{Z})$.    
\end{enumerate}
\ethm	
\begin{proof}
	(1) Let $V \ncong L(X)^A$ for any irreducible uniformly bounded $E^0$-module $X$ and $A \in GL(2, \mathbb{Z})$.
	Using Lemma \ref{L10.14}, pick $\underline{c} = (c_1,c_2) \in
	\mathbb{Z}^2(+)$ with $\Lambda_0 + \underline{c} \in P_D(V)$.
	Now for any $k \in \mathbb{N}$, there
	exists $(a,b) \in \mathbb{Z}^2$ satisfying $0 < a\alpha + b <
	\dfrac{c_1\alpha + c_2}{4k}$, since $>_{\alpha}$ is a dense order
	on $\mathbb{Z}^2$ and $\underline{c} \in \mathbb{Z}^2(+)$. Due to (R3), it is clear that there exist
	\begin{align}\label{(7.14)}
        \text{infinitely many} \ a \in \mathbb{Z} \ \text{and infinitely many} \ b \in \mathbb{Z} \ \text{with} \ (a,b) \ \text{satisfying the above relation}.
	\end{align} 
	Putting $\underline{m} = (a,b) \in \mathbb{Z}^2(+)$, set   
\begin{align*}
l = \max \{t \in \mathbb{Z} \ | \ \Lambda_0 + t \underline{c} \in P_D(V)
\}, \ l^{\prime} = \max \{t \in \mathbb{Z} \ | \ \Lambda_0 +
l\underline{c} + t \underline{m} \in P_D(V) \}. 
\end{align*}
	As $\Lambda_0 + \underline{c} \in P_D(V)$, it follows from (\ref{(4.10)}) that $l \in \mathbb{N}$ and $l^{\prime} \in \mathbb{Z}_{+}$. Let $\mu^{\prime} = \Lambda_0 + l\underline{c} + l^{\prime}\underline{m} = \lambda + \delta_{\underline{s}^{\prime}}$ for some $\underline{s}^{\prime} \in \mathbb{Z}^2$ and $(a_i^{\prime},b_i^{\prime}) = (1- 2i)\underline{m} \ \forall \ 1 \leqslant i \leqslant k$. We can now check that $-\underline{c} <_{\alpha} (a_i^{\prime},b_i^{\prime}) <_{\alpha} (0,0)$. Furthermore set $\underline{p_i}^{\prime} = (2i-1)\underline{m} >_{\alpha} (0,0)$ and  $\underline{q_i}^{\prime} = \underline{c} - \underline{p_i}^{\prime} >_{\alpha} (0,0)$ for each $1 \leqslant i \leqslant k$. Moreover as $V_{\mu^{\prime}}$ is a finite-dimensional $\mathfrak{g}$-module, $V_{\mu^{\prime}}$ has a minimal weight $\eta^{\prime} \in P_{\mathfrak{h}}(V_{\mu^{\prime}})$. Then $\eta^{\prime}(\gamma^{\vee}) \in \{0,1 \}$ for any $\gamma \in \Delta_{\mathfrak{g}}^{+}$, where $\text{span} \{x_{\gamma}, y_{\gamma}, \gamma^{\vee} \} \cong \mathfrak{sl}_2$  Fix any $0 \neq v_0 \in V_{\eta^{\prime} + \delta_{\underline{s}^{\prime}}}$. 
	
	\noindent \textbf{Claim 1.} $K_0$ and $K_1$ act trivially on $V_{\mu^{\prime}}$.\\
	By Lemma \ref{L7.4}, it suffices to only consider $K_0$. If not, then $K_0$ acts by some non-zero $c_0 \in \mathbb{C}$. Now for every $k \in \mathbb{N}$, $\{(x_{\gamma} \otimes t^{-\underline{q_i}^{\prime}})(y_{\gamma} \otimes t^{-\underline{p_i}^{\prime}})v_0 \}_{1 \leqslant i \leqslant k} \subseteq V_{\mu^{\prime} - \underline{c}}$.
	Next we show that the collection $\{(x_{\gamma} \otimes t^{-\underline{q_i}^{\prime}})(y_{\gamma} \otimes t^{-\underline{p_i}^{\prime}})v_0 \}_{1 \leqslant i \leqslant k}$ is not linearly independent for each $k \in \mathbb{N}$.
	If not, then $\dim V_{\mu^{\prime} - \underline{c}} \geqslant k$ for all $k \in \mathbb{N}$. Clearly $(l-1)\underline{c} + l^{\prime}\underline{m} \in \mathbb{Z}^2(+) \cup \{(0,0) \}$. If $(l-1)\underline{c} + l^{\prime}\underline{m} = \underline{0}$, then the assertion is obvious. On the other hand, if $(l-1)\underline{c} + l^{\prime}\underline{m} \in \mathbb{Z}^2(+)$, then we can apply Lemma \ref{L10.13} to infer that, for each $k \in \mathbb{N}$, there exists an injective linear map
	\begin{align*}
	\phi_k:  V_{\mu^{\prime} - \underline{c}} \longrightarrow (V_{\Lambda_0})^{\oplus N}, \ \text{where} \ N= \dim L(\mathfrak{g})_{(1-l)\underline{c} - l^{\prime}\underline{m}}, \ L(\mathfrak{g}) = \mathfrak{g} \otimes \mathbb{C}[t_0^{\pm 1}, t_1^{\pm 1}].
	\end{align*}
	This implies that $N(\dim V_{\Lambda_0}) \geqslant k$ for every $k \in \mathbb{N}$, whence it follows that $\dim V_{\Lambda_0} = \infty$, which is a contradiction. This suggests that we can find some $k^{*} \in \mathbb{N}$ and $ \underline{0} \neq \underline{\beta} = (\beta_1, \ldots, \beta_{k^*}) \in \mathbb{C}^{k^*}$ satisfying
	$\sum_{i=1}^{k^{*}}\beta_i(x_{\gamma} \otimes t^{-\underline{q_i}^{\prime}})(y_{\gamma} \otimes t^{-\underline{p_i}^{\prime}})v_0 = 0$. Now setting $\underline{r}^{\prime} = \underline{c} - 2\underline{m}  >_{\alpha} \underline{0}$, note that      
	\begin{align}\label{(4.11)}
	[E_c, \mathcal{C}] = (0), \ (\mu^{\prime} + \underline{r}^{\prime} -\underline{p_j}^{\prime}) - (\mu^{\prime} + \underline{m}) \in \mathbb{Z}_2(+), \ (\mu^{\prime} + \underline{m}) \notin P_D(V), \\
	 \label{(4.12)} (\mu^{\prime} + \underline{r}^{\prime} -\underline{q_i}^{\prime}) - (\mu^{\prime} + \underline{m}) \in \mathbb{Z}_2(+) \cup \{(0,0)\}  \ \forall \ 2 \leqslant i \leqslant k^*, \ 1 \leqslant j \leqslant k^*, \\ \label{(4.13)}
		(\mu^{\prime} + \underline{q_1}^{\prime} -\underline{r}^{\prime}) \notin P_D(V), \ (\mu^{\prime} + \underline{p_1}^{\prime}) \notin P_D(V), \ (\mu^{\prime} + \underline{c} -\underline{r}^{\prime}) \notin P_D(V).
	\end{align}

\noindent An application of $\gamma^{\vee} \otimes t^{\underline{r}^{\prime}}$ to the last equation together with (\ref{(4.11)}), (\ref{(4.12)}), (\ref{(4.13)}) and (R4) gives  
\begin{align} \label{E6}
\beta_1[2(x_{\gamma} \otimes t^{\underline{r}^{\prime} - \underline{q_1}^{\prime}})(y_{\gamma} \otimes t^{- \underline{p_1}^{\prime}}) + a_1^{\prime} \big (\sum_{i=0}^{1} t^{\underline{r}^{\prime} - \underline{q_1}^{\prime}}K_i \big )(y_{\gamma} \otimes t^{- \underline{p_1}^{\prime}})]v_0 \\ +  \beta_1^{\prime}[2\gamma^{\vee} \otimes t^{\underline{r}^{\prime} - \underline{c}} + a_2^{\prime} (t^{\underline{r}^{\prime} - \underline{c}}K_0) +
a_3^{\prime} (t^{\underline{r}^{\prime} - \underline{c}}K_0)]v_0 = 0, 
\end{align}
where $\beta_1^{\prime} = \sum_{i=2}^{k^*} \beta_i$ and some $a_1^{\prime}, a_2^{\prime}, a_3^{\prime} \in \mathbb{C}$. Consequently applying $\gamma^{\vee} \otimes t^{\underline{c} - \underline{r}^{\prime}}$ to the above equation and using (\ref{(4.11)}), (\ref{(4.12)}) and (\ref{(4.13)}), we obtain
\begin{align} \label{E7}
\big (\eta^{\prime}(\gamma^{\vee}) + ac_0(x_{\gamma}|y_{\gamma}) \big) \beta_1 + \big (ac_0(\gamma^{\vee}|\gamma^{\vee}) \big ) \beta_1^{\prime} = 0. 
\end{align}
Again apply $(y_{\gamma} \otimes t^{\underline{p_1}})(x_{\gamma} \otimes t^{\underline{r} - \underline{q_1}})$ to (\ref{E6}), (10.23) and use (\ref{(4.11)}), (\ref{(4.12)}), (\ref{(4.13)}) along with $(x_{\gamma}|x_{\gamma}) = 0 = (y_{\gamma}|y_{\gamma})$ to get
\begin{align} \label{E8}
\big(2 \eta^{\prime}(\gamma^{\vee}) + 2ac_0(x_{\gamma}|y_{\gamma}) - \eta^{\prime}(\gamma^{\vee})^2 + a^2c_0^2(x_{\gamma}|y_{\gamma})^2 \big) \beta_1 + \big(2 \eta^{\prime}(\gamma^{\vee}) + 2ac_0(x_{\gamma}|y_{\gamma}) \big) \beta_1^{\prime} = 0. 
\end{align}
Now considering the cases $\eta^{\prime}(\gamma^{\vee}) = 1$ and $\eta^{\prime}(\gamma^{\vee})=0$ separately and using (\ref{E7}), (\ref{E8}) along with (\ref{(7.14)}), we can deduce (like in Claim 1 of Lemma \ref{L10.14}) that
\begin{align} \label{E9}
\big (1 + ac_0 \big (x_{\gamma}|y_{\gamma}) \big ) \beta_1 + \big (ac_0(\gamma^{\vee}|\gamma^{\vee}) \big )\beta_1^{\prime} = 0 = \big (1 + ac_0(x_{\gamma}|y_{\gamma}) \big)^2 \beta_1 + \big (2 + 2ac_0(x_{\gamma}|y_{\gamma}) \big ) \beta_1^{\prime}. \\
\big (2 + ac_0(x_{\gamma}|y_{\gamma}) \big ) \beta_1 + 2 \beta_1^{\prime} = 0 = (x_{\gamma} | y_{\gamma}) \beta_1 + (\gamma^{\vee} | \gamma^{\vee}) \beta_1^{\prime}. 
\end{align}
Finally since $c_0(x_{\gamma}|y_{\gamma})(\gamma^{\vee}|\gamma^{\vee}) \neq 0$, we can invoke (\ref{(7.14)}) to pick some $a \in \mathbb{Z}$ such that the system of linear equations in both (\ref{E9}) and (10.27) has a unique solution. This implies that $\beta_1 = \beta_1^{\prime} = 0$. We can similarly deduce that $\beta_2 = \ldots = \beta_{k^*} = 0$, which is a contradiction. Hence the claim. 	
	
\noindent \textbf{Claim 2.} $\mathfrak{h}$ acts trivially on $V_{\mu^{\prime}}$. \\
If not, then there exists $0 \neq v_0 \in V_{\mu^{\prime}}$ with $\gamma^{\vee}v_0 = \mu^{\prime}(\gamma^{\vee})v_0 \neq 0$ for some $\gamma^{\vee} \in \mathfrak{h}, \ \gamma \in \Delta_{\mathfrak{g}}^{+}$. As in the previous claim, there exist some $k^* \in \mathbb{N}$ and $\underline{0} \neq (\alpha_1, \ldots, \alpha_{k^*}) \in \mathbb{C}^{k^*}$ such that $\sum_{i=1}^{k^*}\alpha_i(x_{\gamma} \otimes t^{-\underline{q_i}^{\prime}})(y_{\gamma} \otimes t^{-\underline{p_i}^{\prime}})v_0 = 0$. Applying $(\gamma^{\vee} \otimes t^{\underline{c} - \underline{r}^{\prime}})(\gamma^{\vee} \otimes t^{\underline{r}^{\prime}})$ to this equation and using (\ref{(4.11)}), (\ref{(4.12)}), (\ref{(4.13)}), it can be deduced using Claim 1 (see (\ref{E7})) that $\alpha_1 \mu^{\prime}(\gamma^{\vee}) = 0$, which implies that $\alpha_1 = 0$. It can be similarly shown that $\alpha_2 = \ldots = \alpha_{k^*} = 0$, which is a contradiction. Hence the claim. 

\noindent \textbf{Claim 3.} $(x_{\gamma} \otimes t^{-\underline{m}})v = 0 = (y_{\gamma} \otimes t^{-\underline{m}})v \ \forall \ \gamma \in \Delta_{\mathfrak{g}}^{+}, \ v \in V_{\mu^{\prime}}$. \\
Consider any $\gamma \in \Delta_{\mathfrak{g}}^{+}$ and fix $0 \neq v_0 \in V_{\mu^{\prime}}$. As in Claim 1, there exists $k^* \in \mathbb{N}$ and $\underline{0} \neq \underline{\nu} = (\nu_1, \ldots, \nu_{k^*}) \in \mathbb{C}^{k^*}$ such that $\sum_{i=1}^{k^*}\nu_i(x_{\gamma} \otimes t^{-\underline{q_i}^{\prime}})(x_{\gamma} \otimes t^{-\underline{p_i}^{\prime}})v_0 = 0$. Without loss of generality, take $\nu_1 \neq 0$. Applying $(\gamma^{\vee} \otimes t^{\underline{r}^{\prime}})$ to this equation and using (\ref{(4.11)}), (\ref{(4.12)}), (\ref{(4.13)}), we get
\begin{align*}
\nu_1 \big (2(x_{\gamma} \otimes t^{\underline{r}^{\prime} - \underline{q_1}^{\prime}})(x_{\gamma} \otimes t^{-\underline{m}}) + \sum_{i=0}^{1} \eta_i(t^{\underline{r}^{\prime}- \underline{q_1}^{\prime}}K_i)(x_{\gamma} \otimes t^{-\underline{m}}) \big )v_0 + \nu_1^{\prime} \big (\sum_{i=0}^{1} \eta_i^{\prime}(t^{\underline{r}^{\prime} - \underline{c}}K_i) \big )v_0 = 0,
\end{align*} 
where $\nu_1^{\prime} = \sum_{i=2}^{k^*} \nu_i$ and $\eta_0, \eta_0^{\prime}, \eta_1, \eta_1^{\prime} \in \mathbb{C}$. A further application of $(y_{\gamma} \otimes t^{\underline{m}})$ to the above equation along with (\ref{(4.11)}), Claim 1 and Claim 2 gives us $(x_{\gamma} \otimes t^{-\underline{m}})v_0 = 0$. We can similarly obtain $(y_{\gamma} \otimes t^{-\underline{m}})v_0 = 0$, by simply interchanging the roles of $x_{\gamma}$ and $y_{\gamma}$. Hence the claim. 

\noindent \textbf{Claim 4.} $(\mathfrak{g} \otimes \mathbb{C}t^{-2\underline{m}})v_0 = (0)$. \\
By Lemma \ref{L3.6}, we have $a(t^{-l\underline{m}}K_0) + b (t^{-l\underline{m}}K_1) = 0$ for any $l \in \mathbb{N}$, whence $(\gamma^{\vee} \otimes t^{-2\underline{m}})v_0 = 0$ by Claim $3$. Consequently  $0 =(x_{\gamma} \otimes t^{-\underline{m}})(\gamma^{\vee} \otimes t^{-2\underline{m}})v_0 = (-2x_{\gamma} \otimes t^{-3\underline{m}})v_0$, which implies that $(x_{\gamma} \otimes t^{-3\underline{m}})v_0 = 0$. Again considering $y_{\gamma}$ instead of $x_{\gamma}$ yields $(y_{\gamma} \otimes t^{-3\underline{m}})v_0 = 0$. Then it follows from (\ref{(4.11)}) that $(2x_{\gamma} \otimes t^{-2\underline{m}})v_0 = [\gamma^{\vee} \otimes t^{\underline{m}}, x_{\gamma} \otimes t^{-3 \underline{m}}]v_0 = [y_{\gamma} \otimes t^{-3 \underline{m}}, \gamma^{\vee} \otimes t^{\underline{m}}]v_0 = (2y_{\gamma} \otimes t^{-2\underline{m}})v_0 = 0$. As a result, $(\mathfrak{g} \otimes \mathbb{C}t^{-2\underline{m}})v_0 = (0)$, which proves the claim. 

But as $2\underline{m} \in \mathbb{Z}^2(+)$, this clearly contradicts Lemma \ref{L10.13} and hence we get the desired result. \\
(2) From (1), we already know that $V \cong L(X)^A$ for a uniformly bounded irreducible module $X$ over $E^0$ and $A \in GL(2,\mathbb{Z})$. Again by Claim 2 of Lemma \ref{L8.3}, $X$ is an integrable level zero module over $\mathfrak{g} \otimes \mathbb{C}[t_1,t_1^{-1}] \oplus \mathbb{C}K_1 \oplus \mathbb{C}d_1$. Applying Lemma \ref{L9.1}, we thereby get $0 \neq x \in X$ such that $(\mathfrak{n}^+ \otimes \mathbb{C}[t_1,t_1^{-1}])x = (0)$ and consequently $(E_{+})x = (0)$. Finally using the PBW theorem and the irreducibility of the quasi-finite module $V$, we can deduce that $(0) \neq T := \{v \in V \ | \ (E_{+})v = (0) \}$ is an irreducible level zero quasi-finite module over $E_0$. This implies that $V \cong V(T)^B$ for some $B \in GL(2, \mathbb{Z})$. The theorem now directly follows from Lemma \ref{Lemma} and Remark \ref{Notation}.   
\end{proof}

\brmk
$V(\psi,a)$ need not be a quasi-finite module over $E$ for any $\psi \in (E_0^{\prime})^*$ satisfying $\psi(K_1) = 0$. The precise description of the elements $\psi \in (E_0^{\prime})^*$ for which $V(\psi,a)$ is quasi-finite depends heavily on the Lie algebra $E$. We illustrate this fact by means of two examples. If $E = \tau_{2}^{M}$ (see \S \ref{SS4.4}) and $V(\psi,a)$ is a quasi-finite $E$-module, then it follows from \cite[Corollary 1.24]{E3} that $\psi(K_0)= 0$. On the other hand, if $E = \tau_2^{\mathcal{S}}$ (see \S \ref{SS4.3}), then \cite[Theorem 5.5]{YB} reveals that there do exist quasi-finite modules $V(\psi,a)$ over $E$ with $\psi(K_0) \neq 0$. 
Also see \cite{E2} for examples of irreducible non-zero level quasi-finite modules of highest weight type over $\tau_2^{tor}$.       
\ermk



\section{Main results} \label{S11}
In this section, we use results from Section \ref{S8}, Section \ref{S9} and Section \ref{S10} to classify the irreducible quasi-finite (equivalently weakly integrable by Lemma \ref{L7.8}) modules with non-trivial core action over the toroidal Lie algebra $\tau_2^{tor}$ and fgc EALAs of nullities $1$ and $2$. We also broadly describe these irreducible quasi-finite modules over $\tau_{n+1}^{tor}$ ($n>1$) and fgc EALAs of nullity greater than $2$.

\bthm \label{T11.1}  
Let $V$ be an irreducible quasi-finite module over $\widehat{L}(\mathfrak{g})$. Then $V$ is either isomorphic to an irreducible integrable loop module or a highest weight module $V(\Lambda)$ or its graded dual $V(\Lambda)^*$, for some $\Lambda \in H^*$ with $\Lambda|_{\mathfrak{h}} \in P_{\mathfrak{g}}^{+}$.  
\ethm  

\begin{proof}
The result follows directly from Theorem \ref{T8.5}.
\end{proof}
\brmk \label{General}
\
\begin{enumerate}
	\item The above theorem completely classifies the irreducible \textit{weakly integrable} $\widehat{L}(\mathfrak{g})$-modules, thereby generalizing the results of Chari \cite{C} and Chari--Pressley \cite{CP}, where they have classified the irreducible \textit{integrable} modules over $\widehat{L}(\mathfrak{g})$. 
	\item In \cite[Proposition 7.1(i)]{EF}, the irreducible \textit{weakly integrable} modules of \textit{non-zero level} were shown to be highest weight modules for $\widehat{L}(\mathfrak{g})$. Our Theorem \ref{T11.1} generalizes this result, due to which we finally obtain the complete classification of irreducible weakly integrable $\widehat{L}(\mathfrak{g})$-modules in case of \textit{both level zero and non-zero level} modules.
\end{enumerate}
\ermk

\bthm \label{T11.3} 
Let $V$ be an irreducible quasi-finite module over $E:= E_2$ (or over $\tau_2^{tor}$) with non-trivial action of $E_c$ (respectively $\overline{L}(\mathfrak{g})$). Then $V$ satisfies the following.
\begin{enumerate}
	\item $V$ is either a level zero integrable module or a highest weight type module.
	\item If $V$ is a level zero integrable module, then $\mathcal{C}_2$ (respectively $\mathcal{Z}_2$) acts trivally on $V$ and $V$ is a highest weight module induced from an irreducible uniformly bounded module over $\mathfrak{h} \otimes \mathbb{C}[t_0^{\pm 1}, t_1^{\pm 1}] \oplus \mathcal{D}$ (respectively over $\mathfrak{h} \otimes \mathbb{C}[t_0^{\pm 1}, t_1^{\pm 1}] \oplus D$).  
	\item A highest weight type module is isomorphic to $V(\psi,a)^B$ for some $\psi \in (E_0^{\prime})^*$ (respectively $((\tau_2^{tor})_0^{\prime})^*$) with $\psi(K_1) = 0, \ a \in \mathbb{C}$ and $B \in GL(2, \mathbb{Z})$.
   \end{enumerate}
\ethm

\begin{proof}
(1) follows from Proposition \ref{P9.9} and Theorem \ref{T10.15}, whereas (2) and (3) are immediate consequences of Proposition \ref{P9.5} and Theorem \ref{T10.15} respectively.
\end{proof}

\brmk \label{R11.4}
\
\begin{enumerate}
\item If $V$ is an irreducible quasi-finite module over $E:=E_2$ (or over $\tau_2^{tor}$) of \textit{non-zero level}, then it follows from Theorem \ref{T11.3} that $V \cong V(\psi,a)^B$ for some $\psi \in (E_0^{\prime})^*$ (respectively $((\tau_2^{tor})_0^{\prime})^*$) with $\psi(K_0) \neq 0, \ \psi(K_1) = 0, \ a \in \mathbb{C}$ and $B \in GL(2, \mathbb{Z})$. 
\item The \textit{level zero} irreducible integrable modules over $\tau_{2}^{tor}, \ \tau_{2}^{\mathcal{S}}$ and $\tau_{2}^{M}$ (see Section \ref{S4}) with non-trivial core action were classified in \cite{E1}, \cite{CLT} and \cite{E1} respectively, which in turn gives us concrete realizations of the highest weight spaces of the level zero irreducible integrable modules (occuring in Theorem \ref{T11.3}) over all these Lie algebras (see Remark \ref{R9.6}).
\end{enumerate}       
\ermk

\bthm \label{T11.5}
Let $V$ be an irreducible quasi-finite module over $E_{n+1}$ (or over $\tau_{n+1}^{tor}$) with non-trivial action of $(E_{n+1})_c$ (respectively $\overline{L}(\mathfrak{g})$), where $n>1$. Then $V$ satisfies the following.
\begin{enumerate}
\item $V$ is either a level zero integrable module or a restricted generalized highest weight module.
\item If $V$ is a level zero integrable module, then $\mathcal{C}_{n+1}$ (respectively $\mathcal{Z}_{n+1}$) acts trivally on $V$ and $V$ is a highest weight module induced from an irreducible uniformly bounded module over $\mathfrak{h} \otimes \mathbb{C}[t_0^{\pm 1}, \ldots, t_n^{\pm 1}] \oplus \mathcal{D}$ (respectively over $\mathfrak{h} \otimes \mathbb{C}[t_0^{\pm 1}, \ldots, t_n^{\pm 1}] \oplus D$).  
\item If $V$ is a restricted generalized highest weight (GHW) module with a restricted GHW vector $v_0$ having weight $\Lambda_0$, then $V$ satisfies the properties as stated in Lemma \ref{L9.10}. 
\end{enumerate}
\ethm

\begin{proof}
(1) follows from Proposition \ref{P9.9}, while (2) and (3) are applications of Proposition \ref{P9.5} and Lemma \ref{L9.10} respectively.
\end{proof}
\brmk \label{R11.6}
\
\begin{enumerate}
\item The level zero irreducible integrable modules over $\tau_{n+1}^{tor}, \ \tau_{n+1}^{\mathcal{S}}, \ \tau_{n+1}^{\mathcal{H}}$ and $\tau_{n+1}^{M}$ (see Section \ref{S4}) with non-trivial core action were classified in \cite{E1}, \cite{TB}, \cite{E4} and \cite{E1} respectively. Combining these results along with Theorem \ref{T11.5}, we obtain a description of the irreducible quasi-finite modules with non-trivial core action over all these Lie algebras.    
\item If $V$ is an irreducible quasi-finite module over the toroidal EALA $\tau_{n+1}^{\mathcal{S}}$ (see \S \ref{SS4.3}) or over the Hamiltonian EALA  $\tau_{n+1}^{\mathcal{H}}$ (see \S \ref{SS4.5}) with the \textit{trivial} core action, then $V$ turns out to be irreducible over $\mathcal{S}_{n+1}$ (respectively $\mathcal{H}_{n+1}$). The problem of classifying irreducible quasi-finite modules over $\mathcal{S}_{n+1}$ and $\mathcal{H}_{n+1}$ have been studied by many authors \cite{BT,GL,JL,LT,T}. On the other hand, if $V$ is irreducible over $\tau_{n+1}^{M}$ (or $\tau_{n+1}^{tor}$) with the \textit{trivial} action of the core (respectively $\overline{L}(\mathfrak{g})$), then it is clear that $V$ is isomorphic to a $1$-dimensional highest weight module $V(\delta_{\underline{k}})$ for some $\underline{k} \in \mathbb{Z}^{n+1}$.           
\end{enumerate}
\ermk

\smallskip

\noindent \textbf{Conflict of interest.} The author declares that there is no conflict of interest.

\smallskip

\noindent \textbf{Data availability.} Data sharing is not applicable to this article as no new data were created or analyzed in this study.

\end{document}